\documentclass{article}
\usepackage{amsmath,amssymb,amsthm}
\usepackage{graphicx}
\usepackage[utf8]{inputenc}
\usepackage{anysize}
\usepackage{enumerate}
\marginsize{2cm}{2cm}{2cm}{3cm}

\usepackage{hyperref} 
\hypersetup{colorlinks=true, citecolor=blue, urlcolor=blue} 

\usepackage{tikz}
\usepackage{colortbl}
\usepackage{stmaryrd} 

\theoremstyle{definition}
  \newtheorem{definition}{Definition}[section]
  \newtheorem{remark}{Remark}

\theoremstyle{plain}
  \newtheorem{theorem}{Theorem}[section]
  \newtheorem{proposition}{Proposition}[section]

\newcommand{\re}{\mathbb{R}}

\usepackage{dsfont}
\usepackage{multirow}

\usepackage{caption}
\usepackage{subcaption}

\title{On a class of nonautonomous quasilinear systems  
\\
with general time-gradually-degenerate damping
}
\author{Richard De la cruz\thanks{School of Mathematics and Statistics, Universidad Pedagógica y Tecnológica de Colombia, 150003, Tunja, Colombia. E-mail: \href{richard.delacruz@uptc.edu.co}{richard.delacruz@uptc.edu.co}}  
\,and
Wladimir Neves\thanks{Instituto de Matemática, Universidade Federal do Rio de Janeiro, Cidade Universitária 21945-970, Rio de Janeiro, Brazil. E-mail: \href{wladimir@im.ufrj.br}{wladimir@im.ufrj.br}}
}
\date{\today}

\begin{document}

\maketitle

\begin{abstract}
    In this paper, we study two systems with a time-variable coefficient and general time-gradually-degenerate damping. 
    More explicitly, we construct the Riemann solutions to the time-variable coefficient Zeldovich approximation and 
    time-variable coefficient pressureless gas systems both with general time-gradually-degenerate damping. 
    Applying the method of similar variables and nonlinear viscosity, we obtain classical Riemann solutions and delta shock wave solutions.
\end{abstract}

{\em Keywords:} 
Pressureless gas dynamics system, Zeldovich type approximate system, time-gradually-degenerate damping, Riemann problem,
delta shock solution.

\section{Introduction}

One can find many problems from Continuum Physics that are mathematically modeled by 
balance laws, that is to say, systems of partial differential equations in the following divergence form
\begin{equation}
\label{EQ1}
   \frac{\partial \mathbf{u}}{\partial t} + \sum_{j= 1}^d \frac{\partial \mathbf{F}^j(\mathbf{u})}{\partial x_j}= \mathbf{G}(\mathbf{u}), 
\end{equation}
where $(t,\mathbf{x}) \in \re_+^{d+1} \equiv (0,\infty) \times \re^d$ is the set of independent variables, $\mathbf{u} \in \re^n$ 
denotes the unknown vector field, $\mathbf{F}^j \in \re^n$ is called the flux function and 
$\mathbf{G} \in \re^n$ is the vector production, absorption, or damping term. 
The first component $t>0$ is the 
time variable and $\mathbf{x} \in \re^d$ is the space variable.  
Moreover, when $\mathbf{G} \equiv 0$ equation \eqref{EQ1} 
is called a system of conservation laws. 
In fact, denoting $\mathbf{A}^j(\cdot)= D\mathbf{F}^j(\cdot)$, that is the Jacobian matrices of the fluxes, 
the system \eqref{EQ1} falls in the general class of nonhomogeneous quasilinear first-order systems of partial differential equations
\begin{equation}
\label{EQ2}
   \frac{\partial \mathbf{u}}{\partial t} + \sum_{j= 1}^d \mathbf{A}^j(\mathbf{u}) \frac{\partial \mathbf{u}}{\partial x_j}= \mathbf{G}(\mathbf{u}). 
\end{equation}

Albeit, there are important applications that require to consider systems where the coefficients $\mathbf{A}^j$ and $\mathbf{G}$ in \eqref{EQ2} may depend also on 
the independent variables $(t,\mathbf{x})$, for instance to take into account
material inhomogeneities, or some special geometries, also external actions, etc., see Francesco Oliveri \cite{Francesco Oliveri} and references therein. 
Therefore, one has to study the general nonautonomous quasilinear system of partial differential equations
$$
   \frac{\partial u_i}{\partial t} + \sum_{j= 1}^d A_i^j(t,\mathbf{x}, \mathbf{u}) \frac{\partial \mathbf{u}}{\partial x_j}= G_i(t, \mathbf{x}, \mathbf{u}), \quad (i= 1, \ldots, n). 
$$
We are interested in studying these types of systems, more precisely,
a particular class of such systems which is the $2 \times 2$ systems, $(n=2, d= 1)$, when $A_i \equiv A_i(t,\mathbf{u})$, 
and thus the companion function $G_i= G_i(t,\mathbf{u})$, $(i=1,2)$. Moreover, in this case, we recover in a simple way the 
divergence form. Indeed, taking especially, $A_i(t,\mathbf{u})= \alpha_i(t) A_i(\mathbf{u})$ and $G_i(t,\mathbf{u})= \sigma_i(t) G_i(\mathbf{u})$, 
we may write the above system as 
\begin{equation}
\label{EQQF}
\left\{
\begin{aligned}
\frac{\partial u_1}{\partial t} + \alpha_1(t) \frac{\partial F_1(u_1,u_2)}{\partial x}&= \sigma_1(t) \, G_1(u_1,u_2), 
 \\[5pt]
 \frac{\partial u_2}{\partial t} + \alpha_2(t) \frac{\partial F_2(u_1,u_2)}{\partial x}&= \sigma_2(t) \, G_2(u_1.u_2).
 \end{aligned}
 \right.
\end{equation}

Related to system \eqref{EQQF},
let us start our study by
considering the following
class of nonautonomous quasilinear systems
with time-variable coefficients and time-dependent (linear) damping
represented by the following systems: 
\begin{equation} 
\label{systemLD}
\left \{
\begin{aligned}
    \rho_t + \alpha(t) \, (\rho \, u)_x&= 0,
    \\
    u_t + \alpha(t) \, (\frac{u^2}{2})_x&= - \sigma(t) \, u,
\end{aligned}
\right.
\end{equation}
and also 
\begin{equation} 
\label{ZP}
\left \{
\begin{aligned}
    \rho_t + \alpha(t) \, (\rho \, u)_x&= 0,
    \\
    (\rho u)_t + \alpha(t) \, (\rho u^2)_x&= - \sigma(t) \, \rho u,
\end{aligned}
\right.
\end{equation}
where $0 \le \alpha \in L^1([0,\infty))$, $0 \leq \sigma \in L_{\rm loc}^1([0,\infty))$, the unknown $\rho$ can be interpreted as some density, 
and $u$ is the velocity vector field which carries the density $\rho$.
Companion to \eqref{systemLD} and \eqref{ZP}
the initial data is given by
\begin{equation} \label{DtIni}
    (\rho(x,0),u(x,0))=(\rho_0(x),u_0(x)) =
    \begin{cases}
        (\rho_-,u_-), &\mbox{if } x<0,\\
        (\rho_+,u_+), &\mbox{if } x>0, 
    \end{cases}
\end{equation}
for arbitrary constant states $u_\pm$ and $\rho_\pm >0$.
Therefore, we are considering in fact the Riemann problem, which is the building block of the Cauchy problem. 

\medskip
At this point, we would like to address the reader to \cite{Meng_etal}, where it is studied the following generalized 
Boussinesq system with variable-coefficients, (compare it with the system \eqref{systemLD}), 
$$
\left\{
    \begin{aligned}
        u_t &+ \alpha_1(t) \, (\frac{u^2}{2})_x + \beta_1(t) \,  u_x + \gamma_1(t) \, \rho_x= 0,
        \\[5pt]
        \rho_t & + \alpha_2(t) (\rho \, u)_x + \beta_2(t) \, \rho_x + \gamma_2(t) u_{xxx}= 0,
    \end{aligned}
\right.
$$
where $\alpha_i, \beta_i, \gamma_i$, $(i=1,2)$, are time-dependent coefficients relevant to density,
dispersion and viscosity of the fluid. The above system can model the propagation of weakly dispersive and long 
weakly nonlinear surface waves in shallow water. The authors,
under a selection of the spectral parameters, showed the existence of soliton solutions applying the 
Darboux transformation and symbolic computation.

\medskip
One observes that the second equation of the system \eqref{systemLD} is the Burgers equation with time variable coefficients \cite{DLW}.
In particular,  the time variable coefficients can provide more useful models in many complicated physical situations \cite{DLW, Hong, XuEtal}.
The homogeneous case of the system \eqref{systemLD}, that is to say $\sigma(t)=0$ for all $t \ge 0$, is the following time variable coefficient system
\begin{equation} 
\label{ZeldH}
\left \{
    \begin{aligned}
        \rho_t + \alpha(t) (\rho u)_t&= 0,
        \\
        u_t + \alpha(t) (\frac{u^2}{2})_x&= 0,
    \end{aligned}
    \right.
\end{equation}
which can be interpreted as an extension of {\em Zeldovich approximation system} \cite{VDFN, Zeldovich}. 
In particular, 
the system \eqref{ZeldH} 
with $\alpha(\cdot) \equiv 1$ is used to model the evolution of density inhomogeneities of matter in the universe 
\cite[B. Late nonlinear stage, 3. Sticky dust]{SZ}. Further, let us recall that, the system \eqref{systemLD} belongs to the class of triangular systems 
of conservation laws, that arises in a wide variety of models in physics and engineering, see 
for example \cite{IsTe, SaSch} and references therein. For this reason, the triangular systems have been 
studied by many authors and several rigorous results have been obtained for them. 
In \cite{Delacruz1, Delacruz2}, the Riemann problem was solved to the system \eqref{systemLD} 
with $\alpha(\cdot) \equiv 1$ and $\sigma(\cdot)$ equals to a positive constant, where Delta shocks have to be considered. 
Recently, based on the method of similar variables proposed in \cite{Delacruz2}, Li \cite{Li} studied the Riemann problem to the system \eqref{systemLD}
with $\alpha(\cdot) \equiv 1$ and $\sigma(t)=\frac{\mu}{1+t}$ with physical parameter $\mu >0$.
In the literature the external term $\sigma(t)=\frac{\mu}{(1+t)^\theta}u$ with physical parameters $\mu > 0$ 
and $\theta \ge 0$ is called a {\em time-gradually-degenerate damping} \cite{GengLinMei, LLMZ}, 
and it represents the time-gradually-vanishing friction effect.

\medskip
On the other hand, the homogeneous case of the system \eqref{ZP} is the following time variable coefficient system
\begin{equation} 
\label{ZPH}
\left \{
    \begin{aligned}
        \rho_t + \alpha(t) (\rho \, u)_t&= 0,
        \\
        (\rho u)_t + \alpha(t) (\rho \, u^2)_x&= 0,
    \end{aligned}
    \right.
\end{equation}
which can be seen as an extension of {\em pressureless gas dynamics system} \cite{VDFN, Zeldovich}. 
We recall that gas dynamics with zero pressure is a simplified scenario 
where the pressure of the gas is assumed to be negligible, accounting for
high-speed flows or rarefied gases. 
The first study for the usual pressureless gas dynamics system,
that is \eqref{ZPH} 
with $\alpha(\cdot) \equiv 1$, is
 due to Bouchut \cite{Bouchut} in 1994. In that paper it 
was studied the existence of solutions to the Riemann problem for the pressureless gas dynamics system,
introducing a notion of measure solution and delta shock waves were obtained.
However, uniqueness was not studied.

\medskip
Moreover, the existence of a weak solution to the 
Cauchy problem was first obtained independently by E, Rikov, Sinai \cite{ERS} in 1996, and 
Brenier, Grenier \cite{BG} in 1998. In particular, the authors in \cite{ERS}
show that, the standard entropy condition $(\rho\Phi(\rho))_t + (\rho u\Phi(\rho))_x \le 0$ in the sense of distributions, where $\Phi$ is a convex function, is not enough to express a uniqueness criterion for weak solutions to the Cauchy problem.
Conversely, Wang and Ding \cite{WangDing} proved that the pressureless gas dynamics system has a unique weak 
solution using the Oleinik entropy condition when the initial data $\rho_0$, $u_0$ are both bounded measurable functions.
However, the solution for the Cauchy problem for the pressureless gas dynamics system is in general a Radon measure \cite{ERS}.

\medskip
In 2001, Huang and Wang \cite{HW2001} studied the
Cauchy problem for the system \eqref{ZPH}, when initial data $\rho_0$, $u_0$ are respectively a Radon measure and a bounded measurable function. 
Then, they showed the uniqueness of weak solutions under the Oleinik entropy condition together with an energy condition in the sense that, $\rho u^2$ weakly converges to $\rho_0 u_0^2$ as $t \to 0$.
We recall that, a particular case of Radon measure solution is the delta shock wave solution.
A delta shock wave solution is a type of nonclassical wave solution in which at least one state variable may develop a Dirac measure. 
Actually, on physical grounds, delta shock solutions typically display concentration occurrence in a complex system \cite{Brenier, Korchinski}.
On the other hand, it is well known that the solution for the Riemann problem to the pressureless gas dynamics system involves vacuum and delta shock wave solution and the classical Riemann solutions satisfy the Lax entropy condition while delta shock wave solution is unique under an over-compressive entropy condition \cite{ShenZhang, Yang}.
In a similar way, Keita and Bourgault \cite{Keita} solved the Riemann problem for the pressureless system with linear damping,
that is, the system \eqref{ZP} with $\alpha(\cdot) \equiv 1$ and $\sigma(\cdot) \equiv const.$, showing vacuum states and delta shock solution and uniqueness under the Lax entropy condition and over-compressive entropy condition, respectively.
Finally, De la cruz and Juajibioy \cite{DelacruzJuajibioyAAM2022} obtained delta shock solutions for a generalized pressureless system with linear damping.

\subsection{Equivalent $\times$ non-equivaqlent systems}

Since we considered $\alpha_1= \alpha_2= \alpha$ in both systems \eqref{systemLD},
\eqref{ZP}, one may ask when these systems are equivalent or not. Indeed,  
we observe first that for smooth solutions an elementary manipulation of 
the second equation of \eqref{ZP} reads 
$$ 
    \rho (u_t + \alpha(t) (\frac{1}{2}u^2)_x) + u(\rho_t  +  \alpha(t)(\rho u)_x ) = - \sigma(t) \rho u.
$$
Therefore, due to the first equation of \eqref{ZP} and for $\rho \neq 0$, the above equation reduces to
the equation $(\ref{systemLD})_2$, and thus 
for smooth solutions the system \eqref{ZP} is equivalent to the system \eqref{systemLD}. Albeit, the
question remains open when the solutions are non-regular. 

\medskip
Once placed the above question, we observe that Keita and Bourgault \cite{Keita}, recently in 2019, studied the Riemann problem for the Zeldovich approximation 
and pressureless gas dynamics systems with linear damping with $\sigma = const. >0$. More precisely, 
they analyzed in that paper the Riemann problem to the following systems: 
\begin{equation} 
\label{Kta2}
\left \{
    \begin{aligned}
        \rho_t + (\rho u)_x&= 0,
        \\
        u_t + (\frac{u^2}{2})_x&= - \sigma \, u,
    \end{aligned}
    \right.
\end{equation}
\begin{equation} 
\label{Kta1}
\left \{
    \begin{aligned}
        \rho_t + (\rho u)_x&= 0,
        \\
        (\rho u)_t + (\rho u^2)_x&= - \sigma \rho u,
    \end{aligned}
\right.
\end{equation}
with initial data given by \eqref{DtIni}, and it was proved that
\begin{enumerate}
\item{$u_- < u_+$.} The solution of the Riemann problem \eqref{Kta2}-\eqref{DtIni} and \eqref{Kta1}-\eqref{DtIni} is given by
\begin{equation*}
    (\rho,u)(x,t)=
    \left \{
    \begin{aligned}
        (\rho_-, u_- e^{-\sigma t}), &\qquad x < u_-\frac{1-e^{-\sigma t}}{\sigma},
        \\[5pt]
        (0, \frac{\sigma x}{e^{\sigma t}-1}), &\qquad  u_-\frac{1-e^{-\sigma t}}{\sigma} \le x \le u_+\frac{1-e^{-\sigma t}}{\sigma},
        \\[5pt]
        (\rho_+,u_+ e^{-\sigma t}), &\qquad  x > u_+\frac{1-e^{-\sigma t}}{\sigma}.
    \end{aligned}
    \right. 
\end{equation*}

\item{$u_- > u_+$.}  The solution of the Riemann problem \eqref{Kta2}-\eqref{DtIni} is given by
\begin{equation} 
\label{DeltaSolClassZeld}
    (\rho,u)(x,t)=
    \left \{
    \begin{aligned}
        (\rho_-, u_- e^{-\sigma t}), &\qquad  x < \frac{u_-+u_+}{2\sigma}(1-e^{-\sigma t}),\\
        (w(t) \delta(x-\frac{u_-+u_+}{2\sigma}(1-e^{-\sigma t})), u_\delta(t)), &\qquad x = \frac{u_-+u_+}{2\sigma}(1-e^{-\sigma t}),\\
        (\rho_+,u_+ e^{-\sigma t}), &\qquad x > \frac{u_-+u_+}{2\sigma}(1-e^{-\sigma t}),
    \end{aligned}
    \right. 
\end{equation}
where 
\begin{equation*}
    \label{UDELTA1}
    w(t) = \frac{(\rho_++\rho_-)(u_--u_+)}{2\sigma}(1-e^{-\sigma t}) \quad \text{and} \quad u_\delta(t) = \frac{u_-+u_+}{2}e^{-\sigma t}.
\end{equation*}        
However, the solution of the Riemann problem \eqref{Kta1}-\eqref{DtIni} is given by
\begin{equation*}
    (\rho,u)(x,t)=
    \left \{
    \begin{aligned}
        (\rho_-, u_- e^{-\sigma t}), &\qquad  x < \frac{\sqrt{\rho_+}u_++\sqrt{\rho_-}u_-}{\sqrt{\rho_+}+\sqrt{\rho_-}}(1-e^{\sigma t}),
        \\
        (w(t) \delta(x-\int_0^t u_\delta(s)ds), u_\delta(t)), &\qquad  x = \frac{\sqrt{\rho_+}u_++\sqrt{\rho_-}u_-}{\sqrt{\rho_+}+\sqrt{\rho_-}}(1-e^{\sigma t}),
        \\
        (\rho_+,u_+ e^{-\sigma t}), &\qquad  x > \frac{\sqrt{\rho_+}u_++\sqrt{\rho_-}u_-}{\sqrt{\rho_+}+\sqrt{\rho_-}}(1-e^{\sigma t}),
    \end{aligned}
    \right.
\end{equation*}
where 
\begin{equation*}
    \label{UDELTA2}
w(t) = \frac{\sqrt{\rho_- \rho_+}(u_--u_+)}{\sigma}(1-e^{-\sigma t})
\quad \mbox{ and } \quad
u_\delta(t) = \frac{\sqrt{\rho_+}u_++\sqrt{\rho_-}u_-}{\sqrt{\rho_+}+\sqrt{\rho_-}} e^{-\sigma t}.
\end{equation*} 
\end{enumerate} 

\medskip
Consequently, Keita, Bourgault showed that the systems \eqref{Kta2} and \eqref{Kta1} are equivalent for smooth and also for 
two contact-discontinuity solutions, but they differ for delta shock solutions.
Therefore, it should be expected that a similar scenario of delta shocks are 
presented here as well, and the systems \eqref{systemLD} and
\eqref{ZP} are not equivalent for these types of solutions.  
We remark that the problems here become much more complicated since $\sigma(\cdot)$ besides non-constant 
is just a locally summable function.   

\section{The Zeldovich Type Approximate System} 
\label{Section2}

In this section, we study the Riemann problem to the time-variable coefficient 
Zeldovich's approximate system and time-variable linear damping, that is to say \eqref{systemLD}-\eqref{DtIni}. 
We extended some ideas from \cite{Delacruz1} 
to construct the viscous solutions to the system \eqref{systemLD}, see \eqref{systemVis} below. 
After we show that the family of viscous solutions 
$\{(\rho^\varepsilon,u^\varepsilon)\}$ converges to a solution of the Riemann problem \eqref{systemLD}-\eqref{DtIni}. 
For $u_-<u_+$, classical Riemann solutions are obtained. When $u_->u_+$, we show that a delta shock solution 
is a solution to the Riemann problem \eqref{systemLD}-\eqref{DtIni}.

\subsection{Parabolic regularization}

Given $\varepsilon> 0$, we consider the following parabolic regularization for the system \eqref{systemLD}, 
\begin{equation} 
\label{systemVis}
\left \{
 \begin{aligned}
     \rho_t^\varepsilon + \alpha(t) (\rho^\varepsilon u^\varepsilon)_x&= \varepsilon \beta(t)\rho_{xx}^\varepsilon,
     \\[5pt]
     u_t^\varepsilon + \frac{1}{2} \alpha(t)((u^\varepsilon)^2)_x + \sigma(t)u^\varepsilon&= \varepsilon \beta(t)u_{xx}^\varepsilon, 
 \end{aligned}
 \right. 
\end{equation}
where conveniently we define $\beta(t):= \alpha(t) \exp (-\int_0^t \sigma(s) ds)$. 
We search for $(\rho^\varepsilon,u^\varepsilon)$ be an approximate solution of problem \eqref{systemLD}-\eqref{DtIni}, 
which is defined by the parabolic approximation \eqref{systemVis} with initial data given by
\begin{equation} 
\label{DatRV}
    (\rho^\varepsilon(x,0),u^\varepsilon(x,0)) = (\rho_0(x),u_0(x)), 
\end{equation}
where $(\rho_0,u_0)$ is given by \eqref{DtIni}. 

\medskip
Then, the main issue of this section is to solve problem \eqref{systemVis}  
with initial data \eqref{DatRV}. To this end,
we use the auxiliary function
$u(x,t)=\widehat{u}_x(x,t)e^{-\int_0^t \sigma(\tau)d\tau}$ and a version of Hopf-Cole transformation 
which enable us 
to obtain an explicit solution of the viscous system \eqref{systemVis}-\eqref{DatRV}.
The function $\widehat{u}$ will be explained during the proof of the following 
\begin{proposition} 
\label{Prop2.1}
Under the assumptions on the functions $\alpha, \beta, \sigma$, the explicit solution of the problem \eqref{systemVis}-\eqref{DatRV} is given by 
    $$ 
    \rho^\varepsilon(x,t) = \partial_x W^\varepsilon(x,t) 
    \quad \mbox{ and } \quad
    u^\varepsilon(x,t) = \frac{u_+ b_+^\varepsilon(x,t) + u_- b_-^\varepsilon(x,t)}{b_+^\varepsilon(x,t) + b_-^\varepsilon(x,t)}\exp (-\int_0^t \sigma(s) ds), 
    $$
where
\begin{align*}
    W^\varepsilon(x,t) =& \frac{\rho_- \left( x- u_- \int_0^t \beta(s)ds \right)b_-^\varepsilon(x,t) + \rho_+ \left( x- u_+ \int_0^t \beta(s)ds \right)b_+^\varepsilon(x,t)}{b_-^\varepsilon(x,t)+b_+^\varepsilon(x,t)}
    \\[5pt]
    &+(\rho_+-\rho_-) \frac{(\varepsilon \int_0^t \beta(s)ds)^{1/2} \exp \left( -\frac{x^2}{4\varepsilon \int_0^t \beta(s) ds }\right) }{\pi^{1/2}(b_-^\varepsilon(x,t)+b_+^\varepsilon(x,t))}
\end{align*}
and
$$ b_\pm^\varepsilon(x,t) := \pm \frac{1}{(4\pi \varepsilon \int_0^t \beta(s) ds)^{1/2}} \int_0^{\pm \infty} \exp \left( -\frac{(x-y)^2}{4 \varepsilon \int_0^t \beta(s) ds} - \frac{u_\pm y}{2\varepsilon } \right) dy. $$
\end{proposition}
\begin{proof}
1. Firstly we observe that, if $(\widehat{\rho},\widehat{u})$ solves
\begin{equation} \label{systemVis2}
    \begin{cases}
        \widehat{\rho}_t + \alpha(t) e^{-\int_0^t \sigma(\tau) d\tau} \widehat{\rho}_x \widehat{u}_x = \varepsilon \beta(t) \widehat{\rho}_{xx},\\
        \widehat{u}_t + \frac{1}{2}\alpha(t) e^{-\int_0^t \sigma(\tau) d\tau} (\widehat{u}_x)^2 = \varepsilon \beta(t) \widehat{u}_{xx},
    \end{cases}
\end{equation}
with the initial condition given by
\begin{equation*}
    (\rho(x,0), \widehat{u}(x,0))  =
    \begin{cases}
        (\rho_-x,u_-x), &\mbox{if } x<0,\\
        (\rho_+x,u_+x), &\mbox{if } x>0,
    \end{cases}
\end{equation*}
then $(\rho^\varepsilon,u^\varepsilon)$ defined by $(\widehat{\rho}_x^\varepsilon,\widehat{u}_xe^{-\int_0^t \sigma(\tau) d\tau})$ solves the problem \eqref{systemVis}-\eqref{DatRV}.
Indeed, let us recall the generalized Hopf-Cole transformation, see \cite{Hopf, Delacruz1, Joseph}, that is 
\begin{equation} \label{HCt}
   \begin{cases}
       \widehat{\rho}^\varepsilon = C^\varepsilon e^{\frac{\widehat{u}}{2\varepsilon}},\\
\widehat{u}^\varepsilon = -2\varepsilon \ln (S^\varepsilon).
   \end{cases}  
\end{equation}
Then, from system \eqref{systemVis2} and the generalized Hopf-Cole transformation \eqref{HCt}, we have
\begin{equation} \label{HeatEq}
    \begin{cases}
        C_t^\varepsilon = \varepsilon \beta(t) C_{xx}^\varepsilon,\\
        S_t^\varepsilon = \varepsilon \beta(t) S_{xx}^\varepsilon,
    \end{cases}
\end{equation}
with initial data given by
\begin{equation} \label{IDHeatEq}
    (C^\varepsilon(x,0),S^\varepsilon(x,0)) = 
    \begin{cases}
    (\rho_- x e^{-\frac{u_-x}{2\varepsilon}}, e^{-\frac{u_-x}{2\varepsilon}}), &\mbox{if }x<0, \\
    (\rho_+ x e^{-\frac{u_+x}{2\varepsilon}}, e^{-\frac{u_+x}{2\varepsilon}}), &\mbox{if }x>0.
    \end{cases}
\end{equation}

\medskip
2. Now, the solution to the problem \eqref{HeatEq}-\eqref{IDHeatEq} in terms of the heat kernel is
\begin{equation} \label{EqSolH}
\begin{cases}
    C^\varepsilon(x,t) = a_-^\varepsilon(x,t) + a_+^\varepsilon(x,t),\\
    S^\varepsilon(x,t) = b_-^\varepsilon(x,t) + b_+^\varepsilon(x,t),
\end{cases}
\end{equation}
where
$$
a_\pm^\varepsilon(x,t) :=
\pm \frac{\rho_\pm}{(4\pi \varepsilon \int_0^t \beta(s) ds)^{1/2}} \int_0^{\pm \infty} y \exp \left( -\frac{(x-y)^2}{4 \varepsilon \int_0^t \beta(s) ds} - \frac{u_\pm y}{2\varepsilon} \right) dy
$$
and
$$ b_\pm^\varepsilon(x,t) := \pm \frac{1}{(4\pi \varepsilon \int_0^t \beta(s) ds)^{1/2}} \int_0^{\pm \infty} \exp \left( -\frac{(x-y)^2}{4 \varepsilon \int_0^t \beta(s) ds} - \frac{u_\pm y}{2\varepsilon} \right) dy. $$
Moreover, we have
\begin{equation} \label{Inte1}
    \begin{aligned}
    \int_0^{\pm \infty} \partial_y \left( \exp \left( - \frac{(x-y)^2}{4\varepsilon \int_0^t \beta(s)ds} \right) \right) \exp \left( - \frac{u_\pm y}{2\varepsilon} \right) dy = & - \exp \left( - \frac{x^2}{4\varepsilon \int_0^t \beta(s) ds} \right) \\
    &+ \frac{u_\pm}{2\varepsilon} \int_0^{\pm \infty} \exp \left( - \frac{(x-y)^2}{4\varepsilon \int_0^t \beta(s)ds} - \frac{u_\pm y}{2\varepsilon} \right) dy.
\end{aligned}
\end{equation}
On the other hand, it follows that
\begin{equation} \label{Inte2}
    \begin{aligned}
    \int_0^{\pm \infty} \partial_y \left( \exp \left( - \frac{(x-y)^2}{4\varepsilon \int_0^t \beta(s)ds} \right) \right) \exp \left( - \frac{u_\pm y}{2\varepsilon} \right) dy = & 
    \int_0^{\pm \infty} \frac{(x-y)}{2\varepsilon \int_0^t \beta(s)ds} \exp \left( - \frac{(x-y)^2}{4\varepsilon \int_0^t \beta(s)ds} - \frac{u_\pm y}{2\varepsilon} \right) dy \\
    = &  \frac{x}{2\varepsilon \int_0^t \beta(s)ds} \int_0^{\pm \infty}  \exp \left( - \frac{(x-y)^2}{4\varepsilon \int_0^t \beta(s)ds} - \frac{u_\pm y}{2\varepsilon} \right) dy \\
    &-\int_0^{\pm \infty} \frac{y}{2\varepsilon \int_0^t \beta(s)ds} \exp \left( - \frac{(x-y)^2}{4\varepsilon \int_0^t \beta(s)ds} - \frac{u_\pm y}{2\varepsilon} \right) dy.
\end{aligned}
\end{equation}
Therefore, from \eqref{Inte1} and \eqref{Inte2} we obtain
\begin{equation} \label{EqualA}
    \begin{aligned}
        \int_0^{\pm \infty} y \exp \left( -\frac{(x-y)^2}{4\varepsilon \int_0^t \beta(s) ds} - \frac{u_\pm y}{2\varepsilon} \right) dy =&
        2 \varepsilon \int_0^t \beta(s) ds \cdot \exp \left( - \frac{x^2}{4 \varepsilon \int_0^t \beta(s) ds} \right) \\
        & + \left( x- u_\pm \int_0^t \beta(s) ds \right) \int_0^{\pm \infty} \exp \left( - \frac{(x-y)^2}{4 \varepsilon \int_0^t \beta(s)ds} - \frac{u_\pm y}{2 \varepsilon} \right) dy.
    \end{aligned}
\end{equation}

\medskip
3. Finally, we observe that
\begin{equation*}
    \begin{aligned}
        \partial_x \left( \int_0^{\pm \infty} \exp \left( - \frac{(x-y)^2}{4 \varepsilon \int_0^t \beta(s)ds} - \frac{u_\pm y}{2 \varepsilon} \right) dy \right) =
        - \int_0^{\pm \infty} \partial_y \left( \exp \left( - \frac{(x-y)^2}{4\varepsilon \int_0^t \beta(s) ds} \right) \right) \exp \left( - \frac{u_\pm y}{2\varepsilon} \right) dy
    \end{aligned}
\end{equation*}
and from \eqref{Inte1} we have
\begin{equation} \label{EqualB}
    \begin{aligned}
        \partial_x \left( \int_0^{\pm \infty} \exp \left( - \frac{(x-y)^2}{4 \varepsilon \int_0^t \beta(s)ds} - \frac{u_\pm y}{2 \varepsilon} \right) dy \right)&=
        \exp \left( -\frac{x^2}{4\varepsilon \int_0^t \beta(s)ds} \right) \\
        &- \frac{u_\pm}{2\varepsilon} \int_0^{\pm \infty} \exp \left( -\frac{(x-y)^2}{4\varepsilon \int_0^t \beta(s) ds} - \frac{u_\pm y}{2\varepsilon} \right) dy. 
    \end{aligned}
\end{equation}
Therefore, we may write from \eqref{EqSolH} and \eqref{EqualA} that 
\begin{equation*}
    \begin{aligned}
        C^\varepsilon(x,t)&= \rho_- \left[ - \frac{(\varepsilon \int_0^t \beta(s)ds )^{1/2}}{\pi^{1/2}} \exp \left( -\frac{x^2}{4\varepsilon \int_0^t \beta(s)ds} \right) + \left(x-u_- \int_0^t \beta(s)ds \right) b_-^\varepsilon(x,t;1) \right]
        \\
        &+ \rho_+ \left[ \frac{(\varepsilon \int_0^t \beta(s)ds )^{1/2}}{\pi^{1/2}} \exp \left( -\frac{x^2}{4\varepsilon \int_0^t \beta(s)ds} \right) + \left(x-u_+ \int_0^t \beta(s)ds \right) b_+^\varepsilon(x,t;1) \right].
    \end{aligned}
\end{equation*}
Moreover, from \eqref{EqSolH} and \eqref{EqualB} we have
$$ S_x^\varepsilon(x,t) = -\frac{1}{2\varepsilon} (u_- b_-^\varepsilon(x,t) + u_+ b_+^\varepsilon(x,t)). $$
Applying the generalized Hopf-Cole transformation \eqref{HCt}, it follows that 
$$
\begin{aligned}
    \rho^\varepsilon(x,t)&= \widehat{\rho}_x^\varepsilon(x,t)= (C^\varepsilon(x,t)/S^\varepsilon(x,t))_x,
    \\    
    u^\varepsilon(x,t) &= -2 \varepsilon \frac{S_x^\varepsilon}{S^\varepsilon} \exp (-\int_0^t \sigma(s)ds),
\end{aligned}
$$    
and hence the proof is complete.
\end{proof}

\begin{remark}
One observes that, the solution $(\rho^\varepsilon, u^\varepsilon)$ of  
the problem \eqref{systemVis}-\eqref{DatRV} is absolutely continuous with respect to time $t> 0$, 
and smooth in $x \in \mathbb{R}$. 
\end{remark}

\subsection{The Riemann problem}

In this section, we study the Riemann problem to the system \eqref{systemLD} with $\sigma(t) \geq 0$ for all $t \ge 0$,
which means that the damping can degenerate in some open interval contained in $(0,\infty)$. 

To obtain the Riemann solution 
to the problem \eqref{systemLD} with initial data \eqref{DtIni} we use the viscosity system with time-dependent damping \eqref{systemVis} 
with initial data \eqref{DatRV} and analyze the limit behavior as $\varepsilon \to 0^+$ of the solutions $(\rho^\varepsilon,u^\varepsilon)$ obtained in the previous section. 
To follow, we write $b^\varepsilon_\pm(x,t)$ as 
\begin{align*}
    b_\pm^\varepsilon(x,t) & = \pm \frac{1}{(4\pi \varepsilon \int_0^t \beta(s) ds)^{1/2}} \int_0^{\pm \infty} \exp \left( - \frac{(x-y)^2}{4\varepsilon \int_0^t \beta(s)ds} - \frac{u_\pm y}{2 \varepsilon} \right) dy 
    \\[5pt]
    & = \pm \frac{1}{(\pi B_\varepsilon(t))^{1/2}}
    \exp \left( \frac{-x^2+(x-x_\pm(t))^2}{B_\varepsilon(t)} \right)
    \int_0^{\pm \infty} \exp \left( -\frac{(y+x_\pm(t)-x)^2}{B_\varepsilon(t)} \right) dy 
    \\[5pt]
&= \frac{1}{\pi^{1/2}} \exp \left( \frac{-x^2+(x-x_\pm(t))^2}{B_\varepsilon(t)} \right) \int_{\pm (B_\varepsilon(t))^{1/2}(x_\pm(t)-x)}^{\infty} \exp (-y^2) dy
\\[5pt]
& = \frac{1}{\pi^{1/2}} \exp \left( \frac{-x^2+(x-x_\pm(t))^2}{B_\varepsilon(t)} \right) I_\pm^{\varepsilon,t}, 
\end{align*}
where $x_\pm(t)=u_\pm \int_0^t \beta(s)ds$, $B_\varepsilon(t)=4\varepsilon \int_0^t \beta(s)ds$, and 
$$
     I_\pm^{\varepsilon,t} = \int_{\pm (B_\varepsilon(t))^{1/2}(x_\pm(t)-x)}^{\infty} \exp (-y^2) dy.
$$     
As $\varepsilon \to 0^+$, due to the asymptotic expansion of the (complementary) error function (see \cite{Dingle}), we have
\begin{equation*}
        I_\pm^{\varepsilon,t} = 
\left \{
        \begin{aligned}
    \sum\limits_{n=0}^\infty \frac{(-1)^n (2n)!}{n!} \left( \frac{(B_\varepsilon(t))^{1/2}}{\pm 2(x_\pm(t)-x)} \right)^{2n+1} 
    \!\!\!\exp \left( - \frac{(x_\pm(t)-x)^2}{B_\varepsilon(t)} \right), &\quad \mbox{if } \pm (x_\pm(t)-x) > (B_\varepsilon(t))^{1/2},
    \\
            \frac{1}{2} \pi^{1/2}\quad , &\quad \mbox{if } x_\pm(t) = x,
            \\
            \pi^{1/2} - \sum\limits_{n=0}^\infty \frac{(-1)^n (2n)!}{n!} \left( \frac{(B_\varepsilon(t))^{1/2}}{\mp 2(x_\pm(t)-x)} \right)^{2n+1} 
            \!\!\!\exp \left( - \frac{(x_\pm(t)-x)^2}{B_\varepsilon(t)} \right), &\quad \mbox{if } \pm (x_\pm(t)-x) < -(B_\varepsilon(t))^{1/2},
        \end{aligned}
        \right.
\end{equation*}
and therefore we obtain
\begin{equation} 
\label{approxb}
    b_\pm^\varepsilon(x,t) =
    \left \{
    \begin{aligned}
        \pm \frac{Q_\pm}{\pi^{1/2}} \exp \Big( \frac{x^2}{B_\varepsilon(t)} \Big), &\quad \mbox{if } \pm (x_\pm(t)-x) > (B_\varepsilon(t))^{1/2},
        \\
        \frac{1}{2} \exp \Big( -\frac{x^2}{B_\varepsilon(t)} \Big) \qquad, &\quad \mbox{if } x_\pm(t) = x,\\
        \exp \Big( \frac{-x^2+(x_\pm(t)-x)^2}{B_\varepsilon(t)} \Big) \pm \frac{Q_\pm}{\pi^{1/2}} \exp \Big( -\frac{x^2}{B_\varepsilon(t)} \Big), &\quad \mbox{if } \pm (x_\pm(t)-x) < -(B_\varepsilon(t))^{1/2},
    \end{aligned}
    \right. 
\end{equation}
where 
$$
\begin{aligned}
   Q_\pm &= \sum\limits_{n=0}^\infty \frac{(-1)^n (2n)!}{n!} \Big( \frac{(B_\varepsilon(t))^{1/2}}{2(x_\pm(t)-x)} \Big)^{2n+1} 
   \\[5pt]
   &= \varepsilon^{1/2} \, \left(\frac{\big(\int_0^t \beta(s) ds\big)^{1/2}}{x_\pm(t)-x} - 2 \, \varepsilon  \Big( \frac{\big(\int_0^t \beta(s) ds\big)^{1/2}}{x_\pm(t)-x} \Big)^{3} 
   + 12 \, \varepsilon^2  \Big( \frac{\big(\int_0^t \beta(s) ds\big)^{1/2}}{x_\pm(t)-x} \Big)^{5} - \cdots\right) . 
\end{aligned}
$$

\subsubsection{Classical Riemann solutions: $u_- \le u_+$.} 
\label{Sect2.3.1}
In this case, we have the following
\begin{theorem}
    Suppose that $u_- \le u_+$.
    Let $(\rho^\varepsilon,u^\varepsilon)$ be the solution of the viscosity problem \eqref{systemVis}-\eqref{DatRV}. 
    Then, the limit $$\lim\limits_{\varepsilon \to 0^+} (\rho^\varepsilon(x,t),u^\varepsilon(x,t))= (\rho(x,t),u(x,t))$$
    exists in the sense of distributions, and the pair 
    $(\rho(x,t),u(x,t))$ solves 
    the time-variable coefficient Zeldovich approximate system
    and time-dependent damping \eqref{systemLD} with initial data \eqref{DtIni}.
    In addition, if $u_- < u_+$, then 
    \begin{equation*}
        (\rho(x,t),u(x,t)) = 
        \begin{cases}
            (\rho_-,u_- \exp(-\int_0^t \sigma(s)ds)), &\mbox{if } x<x_-(t),\\
            (0, \frac{x}{\int_0^t \beta(s)ds} \exp(-\int_0^t \sigma(s)ds) ), &\mbox{if } x_-(t)<x<x_+(t),\\
            (\rho_+,u_+ \exp(-\int_0^t \sigma(s)ds)), &\mbox{if } x>x_+(t),
        \end{cases}
    \end{equation*}
    and when $u_-=u_+$, then
    \begin{equation*}
        (\rho(x,t),u(x,t)) =
        \begin{cases}
            (\rho_-,u_- \exp(-\int_0^t \sigma(s)ds)), &\mbox{if } x < x_-(t),\\
        (\rho_+,u_- \exp(-\int_0^t \sigma(s)ds)), &\mbox{if } x > x_-(t).
        \end{cases}
    \end{equation*}
\end{theorem}
\begin{proof}
1. First, let us consider the case $x-x_-(t) < -(B_\varepsilon(t))^{1/2}$. For $\varepsilon> 0$ sufficiently small, due to 
    approximations given by \eqref{approxb}, we may write
    \begin{equation*}
        W^\varepsilon(x,t) \approx \frac{\rho_-(x-x_-(t))c_-^\varepsilon - \frac{\rho_+ (B_\varepsilon(t))^{1/2}}{2\pi^{1/2}}\exp \left(-\frac{x^2}{B_\varepsilon(t)} \right) 
        + (\rho_+-\rho_-)\frac{(B_\varepsilon(t))^{1/2}}{2\pi^{1/2}} \exp \left( - \frac{x^2}{B_\varepsilon(t)} \right) }{c_-^\varepsilon + \frac{(B_\varepsilon(t))^{1/2}}{2\pi^{1/2}(x_+(t)-x)}},
    \end{equation*}
    where $c_-^\varepsilon = \exp \left( \frac{-x^2+(x_-(t)-x)^2}{B_\varepsilon(t)} \right) - \frac{(B_\varepsilon(t))^{1/2}}{2\pi^{1/2}(x_-(t)-x)} \exp \left( - \frac{x^2}{B_\varepsilon(t)} \right)$.
    Therefore, we obtain 
    \begin{equation} \label{approxW1}
        W^\varepsilon(x,t) \approx 
        \frac{\rho_-(x-x_-(t)) \exp \left( \frac{(x_-(t)-x)^2}{B_\varepsilon(t)} \right)}{\exp \left( \frac{(x_-(t)-x)^2}{B_\varepsilon(t)} \right) + \frac{(B_\varepsilon(t))^{1/2}}{2\pi^{1/2}} \left( \frac{1}{x_+(t)-x} - \frac{1}{x_-(t)-x} \right)}
    \end{equation}
    and
    \begin{equation} \label{approxu1}
        u^\varepsilon(x,t) \approx
\frac{\frac{(B_\varepsilon(t))^{1/2}}{2\pi^{1/2}} \left( \frac{u_+}{x_+(t)-x} - \frac{u_-}{x_-(t)-x} \right) + u_- \exp \left( \frac{(x_-(t)-x)^2}{B_\varepsilon(t)} \right)}{\exp \left( \frac{(x_-(t)-x)^2}{B_\varepsilon(t)} \right) + \frac{(B_\varepsilon(t))^{1/2}}{2\pi^{1/2}} \left( \frac{1}{x_+(t)-x} - \frac{1}{x_-(t)-x} \right) } \exp ( -\int_0^t \sigma(s) ds ).
    \end{equation}
    
\medskip
2. Similarly, if $x_-(t)+(B_\varepsilon(t))^{1/2} < x < x_+(t) - (B_\varepsilon(t))^{1/2}$, then we approximate $W^\varepsilon$ as
    \begin{equation*}
        W^\varepsilon(x,t) \approx
        \frac{\rho_-(x-x_-(t))\widehat{c}_-^\varepsilon + \rho_+(x-x_+(t))\widehat{c}_+^\varepsilon + (\rho_+-\rho_-) \frac{(B_\varepsilon(t))^{1/2}}{2\pi^{1/2}} \exp \left( - \frac{x^2}{B_\varepsilon(t)} \right)}{\widehat{c}_-^\varepsilon + \widehat{c}_+^\varepsilon}
    \end{equation*}
    where $\widehat{c}_\pm^\varepsilon = \pm \frac{1}{\pi^{1/2}} \left( \frac{(B_\varepsilon(t))^{1/2}}{2(x_\pm(t)-x)} - \frac{(B_\varepsilon(t))^{3/2}}{4(x_\pm(t)-x)^3} \right) \exp \left( - \frac{x^2}{B_\varepsilon(t)} \right)$. Therefore,
    \begin{equation} \label{approxW2}
        W^\varepsilon(x,t) \approx \frac{B_\varepsilon(t) \left( \frac{\rho_+}{(x_+(t)-x)^2} - \frac{\rho_-}{(x_-(t)-x)^2} \right)}{2 \left( \frac{1}{x_+(t)-x} - \frac{1}{x_-(t)-x} \right) + B_\varepsilon(t) \left( \frac{1}{(x_-(t)-x)^3} - \frac{1}{(x_+(t)-x)^3} \right)}
    \end{equation}
    and
    \begin{equation} \label{approxu2}
        u^\varepsilon(x,t) =
        \frac{
\frac{u_+}{x_+(t)-x} + \frac{u_-}{x-x_-(t)}+\sum\limits_{n=1}^\infty \frac{(-1)^n (2n)! (B_\varepsilon(t))^n}{n! 4^n} \left( \frac{u_+}{(x_+(t)-x)^{2n+1}} + \frac{u_-}{(x-x_-(t))^{2n+1}} \right)
}{
\frac{1}{x_+(t)-x} + \frac{1}{x-x_-(t)}+\sum\limits_{n=1}^\infty \frac{(-1)^n (2n)! (B_\varepsilon(t))^n}{n! 4^n} \left( \frac{1}{(x_+(t)-x)^{2n+1}} + \frac{1}{(x-x_-(t))^{2n+1}} \right)
}
\exp (-\int_0^t \sigma(s)ds ).
    \end{equation}
Moreover, if $x_+(t)-x <-(B_\varepsilon(t))^{1/2}$, then we have
    \begin{equation*}
        W^\varepsilon(x,t) \approx
        \frac{\frac{\rho_-(B_\varepsilon(t))^{1/2}}{2 \pi^{1/2}} \exp \left( - \frac{x^2}{B_\varepsilon(t)} \right) + \rho_+ (x-x_+(t))c_+^\varepsilon + (\rho_+-\rho_-) \frac{(B_\varepsilon(t))^{1/2}}{2\pi^{1/2}} \exp \left( - \frac{x^2}{B_\varepsilon(t)} \right)}{\frac{(B_\varepsilon(t))^{1/2}}{2\pi^{1/2}(x-x_-(t))}\exp \left( - \frac{x^2}{B_\varepsilon(t)} \right) + c_+^\varepsilon}
    \end{equation*}
    where $c_+^\varepsilon = \exp \left( \frac{-x^2+(x_+(t)-x)^2}{B_\varepsilon(t)} \right) - \frac{(B_\varepsilon(t))^2}{2 \pi^{1/2}(x-x_+(t))} \exp \left( - \frac{x^2}{B_\varepsilon(t)} \right)$, and therefore we get
    \begin{equation} \label{approxW3}
        W^\varepsilon(x,t) \approx
        \frac{\rho_+(x-x_+(t)) \exp \left( \frac{(x_+(t)-x)^2}{B_\varepsilon(t)} \right)}{\exp \left( \frac{(x_+(t)-x)^2}{B_\varepsilon(t)} \right) + \frac{(B_\varepsilon(t))^{1/2}}{2\pi^{1/2}} \left( \frac{1}{x_+(t)-x} -  \frac{1}{x_-(t)-x} \right)}
    \end{equation}
    and
    \begin{equation} \label{approxu3}
        u^\varepsilon(x,t) \approx
        \frac{u_+ \exp \left( \frac{(x_+(t)-x)^2}{B_\varepsilon(t)} \right) 
        +
\frac{(B_\varepsilon(t))^{1/2}}{2\pi^{1/2}} \left( \frac{u_+}{x_+(t)-x} - \frac{u_-}{x_-(t)-x} \right)
        }{\exp \left( \frac{(x_+(t)-x)^2}{B_\varepsilon(t)} \right) 
        +
\frac{(B_\varepsilon(t))^{1/2}}{2\pi^{1/2}} \left( \frac{1}{x_+(t)-x} - \frac{1}{x_-(t)-x} \right)} 
\exp (-\int_0^t \sigma(s)ds ).
    \end{equation}
    
\medskip 
3. Now, for the case $u_- <u_+$, from \eqref{approxW1}, \eqref{approxW2}, and \eqref{approxW3} we have
    \begin{equation*}
        \lim_{\varepsilon \to 0+} W^\varepsilon(x,t) = W(x,t)
        =
        \begin{cases}
            \rho_-(x-x_-(t)), &\mbox{if } x<x_-(t), \\
            0, & \mbox{if } x_-(t) < x < x_+(t),\\
            \rho_+(x-x_+(t)), &\mbox{if } x>x_+(t)
        \end{cases}
    \end{equation*}
    and from \eqref{approxu1}, \eqref{approxu2}, and \eqref{approxu3} we have
    \begin{equation*}
        \lim_{\varepsilon \to 0+} u^\varepsilon(x,t) = u(x,t) =
        \begin{cases}
            u_- \exp (-\int_0^t \sigma(s)ds), &\mbox{if } x < x_-(t),\\
            \frac{x}{\int_0^t \beta(s)ds} \exp(-\int_0^t \sigma(s)ds), &\mbox{if } x_-(t) < x < x_+(t),\\
            u_+ \exp (-\int_0^t \sigma(s)ds), &\mbox{if } x > x_+(t).
        \end{cases}
    \end{equation*}
Since $u^\varepsilon(x,t)$
is bounded on compact subsets of $\mathbb{R}^2_+ = \{ (x,t) \, : \, x \in \mathbb{R}, t>0 \}$
    and $u^\varepsilon(x,t) \to u(x,t)$ pointwise as $\varepsilon \to 0+$, then
    $u^\varepsilon(x,t) \to u(x,t)$ in the sense of distribution.
Also, $W^\varepsilon(x,t)$ is bounded on compact subsets of $\mathbb{R}^2_+$ and  $W^\varepsilon(x,t) \to W(x,t)$ pointwise as 
$\varepsilon \to 0+$, then $W^\varepsilon(x,t) \to W(x,t)$ in the sense of distributions and so $W_x^\varepsilon(x,t)$ converges in the distributional 
sense to $W_x(x,t)$. From Proposition \ref{Prop2.1}, we have that $\lim\limits_{\varepsilon \to 0+} \rho^\varepsilon(x,t) = \rho(x,t)$ 
exists in the sense of distribution and
    \begin{equation*}
        \rho(x,t) = W_x(x,t) =
        \begin{cases}
            \rho_-, &\mbox{if } x < x_-(t),\\
            0, & \mbox{if } x_-(t) < x < x_+(t),\\
            \rho_+, &\mbox{if } x > x_+(t).
        \end{cases}
    \end{equation*}
For the case $u_-=u_+$,
we have
\begin{equation*}
    \lim_{\varepsilon \to 0+}(\rho^\varepsilon(x,t),u^\varepsilon(x,t)) = (\rho(x,t),u(x,t)) =
    \begin{cases}
        (\rho_-,u_- \exp(-\int_0^t \sigma(s)ds)), &\mbox{if } x < x_-(t),\\
        (\rho_+,u_- \exp(-\int_0^t \sigma(s)ds)), &\mbox{if } x > x_-(t).
    \end{cases}
\end{equation*}
 Finally, it is not difficult to show that $(\rho(x,t),u(x,t))$ solves \eqref{systemLD}, and thus we omit the details.
\end{proof}

\subsubsection{Delta shock wave solutions: $u_- > u_+$.}

In this section, we study the Riemann problem to the system \eqref{systemLD} with initial data \eqref{DtIni} when $u_- > u_+$. 
Let us recall that, in particular when $\alpha(\cdot)\equiv 1$ and $\sigma(\cdot) \equiv \sigma = const.>0$, 
the solution is not bounded and contains a weighted delta measure supported on a smooth curve (see \cite{Keita}), 
which is a delta shock solution given by \eqref{DeltaSolClassZeld}. 

Here, we have a more general context with similar results. 
Therefore, we first define the meaning of a two-dimensional weighted delta function. 
\begin{definition}
Given $w \in L^1((a,b))$, with $-\infty< a < b< \infty$, 
and a smooth curve $$L \equiv \{ (x(s),t(s)) \, : \, a < s < b \},$$ we say that 
$w(\cdot) \delta_L$ is a two-dimensional weighted delta function supported on $L$, 
when for each test function $\varphi \in C_0^\infty(\mathbb{R} \times [0,\infty))$, 
\begin{equation*}
    \langle w(\cdot) \delta_L, \varphi(\cdot,\cdot) \rangle= \int_a^b w(s) \, \varphi(x(s),t(s)) \, ds. 
\end{equation*}
\end{definition}

\smallskip
Now, the following definition tells us when a pair $(\rho, u)$ 
is a delta shock wave solution to the Riemann problem \eqref{systemLD}-\eqref{DtIni}.
\begin{definition} \label{Def3.2}
    A distribution pair $(\rho,u)$ is called a delta shock wave solution of the problem \eqref{systemLD} and \eqref{DtIni} 
    in the sense of distributions, when there exists a smooth curve $L$ and a function $w(\cdot)$, 
    such that $\rho$ and $u$ are represented respectively by 
    $$
    \rho= \widehat{\rho}(x,t) + w\delta_L, \quad  u= u(x,t)
    $$
with $\widehat{\rho}, u \in L^\infty(\mathbb{R} \times (0,\infty))$, and satisfy
     for each the test function $\varphi \in C_0^\infty(\mathbb{R} \times (0,\infty))$,
    \begin{equation*}
    \left \{
        \begin{aligned}
            < \rho,\varphi_t>+<\alpha \rho u,\varphi_x>=0,
            \\
            \iint_{\mathbb{R}^2_+} \Big( u \varphi_t + \frac{\alpha(t)}{2} u^2 \varphi_x - \sigma(t) u \varphi \Big) \, dx dt= 0,
        \end{aligned}
        \right.
    \end{equation*}
   where
    $$
    <\rho,\varphi>= \iint_{\mathbb{R}^2_+} \widehat{\rho} \, \varphi \, dx dt + \langle w \delta_L, \varphi \rangle, 
    $$
    and
    $$
    <\alpha \rho u,\varphi>= \iint_{\mathbb{R}^2_+} \alpha(t) \, \widehat{\rho} \, u \, \varphi \, dx dt 
    + \langle \alpha(\cdot) w u_\delta \delta_L, \varphi \rangle.
    $$
    Moreover, $u|_L = u_\delta(\cdot)$.
\end{definition}

\bigskip
Placed the previous definitions, we are going to show a solution with a discontinuity on $x= x(t)$ for the system \eqref{systemLD} of the form
\begin{equation*}
    (\rho(x,t),u(x,t))=
    \left \{
    \begin{aligned}
        (\rho_-(x,t),u_-(x,t)), &\quad \mbox{if } x<\gamma(t),\\
        (w(t)\delta_L,u_\delta(t)), &\quad \mbox{if } x=\gamma(t),\\
        (\rho_+(x,t),u_+(x,t)), &\quad \mbox{if } x>\gamma(t),
    \end{aligned}
    \right.
\end{equation*}
where $\rho_\pm(x,t)$, $u_\pm(x,t)$ are piecewise smooth solutions of system \eqref{systemLD}, $\delta_L$ is the Dirac measure 
supported on the curve $\gamma \in C^1$, and $\gamma$, $w$, and $u_\delta$ are to be determined.
Then, we have the following

\begin{theorem}
    Suppose $u_->u_+$. Let $(\rho^\varepsilon,u^\varepsilon)$ be the solution of the problem \eqref{systemVis}-\eqref{DatRV}. Then the limit 
    $$
    \lim_{\varepsilon \to 0+} (\rho^\varepsilon(x,t),u^\varepsilon(x,t)) = (\rho(x,t),u(x,t))
    $$
    exists in the sense of distributions and $(\rho(x,t),u(x,t))$ solves the problem \eqref{systemLD}-\eqref{DtIni}. In addition,
    \begin{equation*}
        (\rho(x,t),u(x,t))=
        \left \{
        \begin{aligned}
            (\rho_-,u_- \exp(-\int_0^t \sigma(s)ds)), &\quad \mbox{if } x<x(t),
            \\
            (w(t)\delta(x-x(t)), \frac{u_-+u_+}{2}\exp(-\int_0^t \sigma(s)ds)), &\quad \mbox{if } x=x(t),\\
            (\rho_+,u_+ \exp(-\int_0^t \sigma(s)ds)), &\quad \mbox{if } x>x(t),
        \end{aligned}
        \right. 
    \end{equation*}
    where 
    $$
    \begin{aligned}
        w(t)&= \frac{1}{2}(\rho_-+\rho_+)(u_--u_+)\int_0^t \alpha(s) \exp (-\int_0^s \sigma(\tau) d\tau) ds, 
        \\
        x(t)&= \frac{u_-+u_+}{2} \int_0^t \alpha(s) \exp (-\int_0^s \sigma(\tau) d\tau) ds.
        \end{aligned}
    $$    
\end{theorem}
\begin{proof}
1. First, since $u_->u_+$, it follows that $x_-(t) > x_+(t)$. For $\varepsilon> 0$ sufficiently small, if $x-x_-(t) > (B_\varepsilon(t))^{1/2}$, then we may write from \eqref{approxb}, 
    \begin{equation*}
        W^\varepsilon(x,t) \approx
        \frac{\rho_+(x-x_+(t)) \exp \left( \frac{(x_+(t)-x)^2}{B_\varepsilon(t)} \right)}{\exp \left( \frac{(x_+(t)-x)^2}{B_\varepsilon(t)} \right) + \frac{(B_\varepsilon(t))^{1/2}}{2\pi^{1/2}} \left( \frac{1}{x_+(t)-x} - \frac{1}{x_-(t)-x} \right)} 
    \end{equation*}
    and 
    \begin{equation*}
        u^\varepsilon(x,t) \approx
        \frac{
        u_+ \exp \left( \frac{(x_+(t)-x)^2}{B_\varepsilon(t)} \right) + \frac{(B_\varepsilon(t))^{1/2}}{2\pi^{1/2}} \left( \frac{u_+}{x_+(t)-x} - \frac{u_-}{x_-(t)-x} \right)
        }{
        \exp \left( \frac{(x_+(t)-x)^2}{B_\varepsilon(t)} \right) + \frac{(B_\varepsilon(t))^{1/2}}{2\pi^{1/2}} \left( \frac{1}{x_+(t)-x} - \frac{1}{x_-(t)-x} \right)
        } 
        \exp (-\int_0^t \sigma(s) ds).
    \end{equation*}
If $x_+(t) -x < -(B_\varepsilon(t))^{1/2}$ and $x=x_-(t)$, then
    \begin{equation*}
        W^\varepsilon(x,t) \approx
        \frac{
        -\rho_- \frac{(B_\varepsilon(t))^{1/2}}{2\pi^{1/2}} \exp \left( -\frac{x^2}{B_\varepsilon(t)} \right) + \rho_+ (x-x_+(t)) \exp \left( \frac{-x^2+(x_+(t)-x)^2}{B_\varepsilon(t)} \right)
        }{
        \frac{1}{2} \exp \left( -\frac{x^2}{B_\varepsilon(t)} \right) + \exp \left( \frac{-x^2+(x_+(t)-x)^2}{B_\varepsilon(t)} \right) + \frac{(B_\varepsilon(t))^{1/2}}{2 \pi^{1/2}(x_+(t)-x)} \exp \left( -\frac{x^2}{B_\varepsilon(t)} \right)
        }
    \end{equation*}
    and
    \begin{equation*}
        u^\varepsilon(x,t) \approx
        \frac{
        \frac{u_-}{2} \exp \left( -\frac{x^2}{B_\varepsilon(t)} \right) + u_+ \exp \left( \frac{-x^2+(x_+(t)-x)^2}{B_\varepsilon(t)} \right) +u_+ \frac{(B_\varepsilon(t))^{1/2}}{2\pi^{1/2}(x_+(t)-x)} \exp \left( -\frac{x^2}{B_\varepsilon(t)} \right)
        }{
        \frac{1}{2} \exp \left( -\frac{x^2}{B_\varepsilon(t)} \right) + \exp \left( \frac{-x^2+(x_+(t)-x)^2}{B_\varepsilon(t)} \right) + \frac{(B_\varepsilon(t))^{1/2}}{2 \pi^{1/2}(x_+(t)-x)} \exp \left( -\frac{x^2}{B_\varepsilon(t)} \right)
        } \exp (-\int_0^t \sigma(s)ds ).
    \end{equation*}
    If $x_+(t) + (B_\varepsilon(t))^{1/2} \le x \le x_-(t) - (B_\varepsilon(t))^{1/2}$, then
    \begin{equation*}
        W^\varepsilon(x,t) \approx
        \frac{
        \rho_- (x-x_-(t)) \exp \left( \frac{(x_-(t)-x)^2}{B_\varepsilon(t)} \right)
        +
        \rho_+ (x-x_+(t)) \exp \left( \frac{(x_+(t)-x)^2}{B_\varepsilon(t)} \right)
        }{
        \exp \left( \frac{(x_-(t)-x)^2}{B_\varepsilon(t)} \right)
        +
        \exp \left( \frac{(x_+(t)-x)^2}{B_\varepsilon(t)} \right)
        +
        \frac{(B_\varepsilon(t))^{1/2}}{2\pi^{1/2}} \left( \frac{1}{x_+(t)-x} - \frac{1}{x_-(t)-x} \right)
        }
    \end{equation*}
    and
    \begin{equation*}
        u^\varepsilon(x,t) \approx
        \frac{
        u_- \exp \left( \frac{(x_-(t)-x)^2}{B_\varepsilon(t)} \right)
        +
        u_+ \exp \left( \frac{(x_+(t)-x)^2}{B_\varepsilon(t)} \right)
        +
        \frac{(B_\varepsilon(t))^{1/2}}{2\pi^{1/2}} \left( \frac{u_+}{x_+(t)-x} - \frac{u_-}{x_-(t)-x} \right)
        }{
        \exp \left( \frac{(x_-(t)-x)^2}{B_\varepsilon(t)} \right)
        +
        \exp \left( \frac{(x_+(t)-x)^2}{B_\varepsilon(t)} \right)
        +
        \frac{(B_\varepsilon(t))^{1/2}}{2\pi^{1/2}} \left( \frac{1}{x_+(t)-x} - \frac{1}{x_-(t)-x} \right)
        } \exp (-\int_0^t \sigma(s)ds ).
    \end{equation*}
    If $x_+(t)-x > (B_\varepsilon(t))^{1/2}$, then
    \begin{equation*}
        W^\varepsilon(x,t) \approx
        \frac{
        \rho_- (x-x_-(t)) \exp \left( \frac{(x_-(t)-x)^2}{B_\varepsilon(t)} \right)
        }{
        \exp \left( \frac{(x_-(t)-x)^2}{B_\varepsilon(t)} \right)
        +
        \frac{(B_\varepsilon(t))^{1/2}}{2\pi^{1/2}} \left( \frac{1}{x_+(t)-x} - \frac{1}{x_-(t)-x} \right)
        }
    \end{equation*}
    and
    \begin{equation*}
        u^\varepsilon(x,t) \approx
        \frac{
        u_- \exp \left( \frac{(x_-(t)-x)^2}{B_\varepsilon(t)} \right)
        +
        \frac{(B_\varepsilon(t))^{1/2}}{2\pi^{1/2}} \left( \frac{u_+}{x_+(t)-x} - \frac{u_-}{x_-(t)-x} \right)
        }{
        \exp \left( \frac{(x_-(t)-x)^2}{B_\varepsilon(t)} \right)
        +
        \frac{(B_\varepsilon(t))^{1/2}}{2\pi^{1/2}} \left( \frac{1}{x_+(t)-x} - \frac{1}{x_-(t)-x} \right)
        } \exp (-\int_0^t \sigma(s)ds ).
    \end{equation*}
    If $x-x_-(t) < -(B_\varepsilon(t))^{1/2}$ and $x=x_+(t)$, then
    \begin{equation*}
        W^\varepsilon(x,t) \approx
        \frac{
        \rho_- (x-x_-(t)) \exp \left( \frac{-x^2+(x_-(t)-x)^2}{B_\varepsilon(t)} \right) + \rho_+ \frac{(B_\varepsilon(t))^{1/2}}{2\pi^{1/2}} \exp \left( -\frac{x^2}{B_\varepsilon(t)} \right)
        }{
        \exp \left( \frac{-x^2+(x_-(t)-x)^2}{B_\varepsilon(t)} \right)
        - 
        \frac{(B_\varepsilon(t))^{1/2}}{2\pi^{1/2}(x_-(t)-x)} \exp \left( -\frac{x^2}{B_\varepsilon(t)} \right)
        +
        \frac{1}{2} \exp \left( -\frac{x^2}{B_\varepsilon(t)} \right)
        }
    \end{equation*}
    and 
    \begin{equation*}
        u^\varepsilon(x,t) \approx
        \frac{
        u_- \exp \left( \frac{-x^2+(x_-(t)-x)^2}{B_\varepsilon(t)} \right) 
        -
        u_- \frac{(B_\varepsilon(t))^{1/2}}{2\pi^{1/2}(x_-(t)-x)} \exp \left( -\frac{x^2}{B_\varepsilon(t)} \right)
        + 
        \frac{u_+}{2} \exp \left( -\frac{x^2}{B_\varepsilon(t)} \right)
        }{
        \exp \left( \frac{-x^2+(x_-(t)-x)^2}{B_\varepsilon(t)} \right)
        - 
        \frac{(B_\varepsilon(t))^{1/2}}{2\pi^{1/2}(x_-(t)-x)} \exp \left( -\frac{x^2}{B_\varepsilon(t)} \right)
        +
        \frac{1}{2} \exp \left( -\frac{x^2}{B_\varepsilon(t)} \right)
        } \exp (-\int_0^t \sigma(s)ds ).
    \end{equation*}
    Therefore, we have that
    \begin{equation*}
        \lim_{\varepsilon \to 0+} W^\varepsilon(x,t) =
        \begin{cases}
            \rho_-(x-x_-(t)), &\mbox{if } (x-x_+(t))^2 - (x-x_-(t))^2 <0,\\
            \rho_+(x-x_+(t)), &\mbox{if } (x-x_+(t))^2 - (x-x_-(t))^2 >0.
        \end{cases}
    \end{equation*}
    Observe that $(x-x_+(t))^2 - (x-x_-(t))^2 = 2(x_-(t)-x_+(t)) (x- \frac{x_-(t)+x_+(t)}{2})$,
    and defining 
    $$\frac{x_-(t)+x_+(t)}{2} =: x(t),$$
we get
\begin{equation*}
        \lim_{\varepsilon \to 0+} W^\varepsilon(x,t) =
        \begin{cases}
            \rho_-(x-x_-(t)), &\mbox{if } x<x(t),\\
            \rho_+(x-x_+(t)), &\mbox{if } x>x(t).
        \end{cases}
    \end{equation*}
    Since $W^\varepsilon(x,t)$ is bounded on compact subsets of $\mathbb{R}^2_+$ and
$W^\varepsilon(x,t) \to W(x,t)$ pointwise as $\varepsilon \to 0+$, then $W^\varepsilon(x,t) \to W(x,t)$ in the sense of distribution and so $W^\varepsilon_x(x,t)$ converges in the distributional sense to $W_x(x,t)$. From Proposition \ref{Prop2.1} we have that $\lim\limits_{\varepsilon \to 0+} \rho^\varepsilon(x,t) = \rho(x,t)$ exists in the sense of distribution and
\begin{equation} \label{Delrho1}
    \rho(x,t) = W_x(x,t) =
    \left \{
    \begin{aligned}
        \rho_- \quad , &\quad \mbox{if } x<x(t),\\
        (x_-(t)-x_+(t)) \frac{\rho_-+\rho_+}{2} \delta(x-x(t)), &\quad \mbox{if } x= x(t)\\
        \rho_+ \quad , &\quad \mbox{if } x>x(t).
    \end{aligned}
    \right. 
\end{equation}
Analogously, we obtain
\begin{equation} \label{Delu1}
    u(x,t)=
    \left \{
    \begin{aligned}
        u_-\exp (-\int_0^t \sigma(s)ds), &\quad \mbox{if } x < x(t),\\
        \frac{u_-+u_+}{2}\exp (-\int_0^t \sigma(s)ds), &\quad \mbox{if } x = x(t),\\
        u_+\exp (-\int_0^t \sigma(s)ds), &\quad \mbox{if } x > x(t).
    \end{aligned}
    \right. 
\end{equation}

\medskip
2. Now, we show that $\rho$ and $u$, defined respectively by \eqref{Delrho1}, \eqref{Delu1} solve the Riemann problem 
\eqref{systemLD}-\eqref{DtIni} in the sense of Definition \ref{Def3.2}. Indeed, for any test function
$\varphi \in C_0^\infty(\mathbb{R} \times \mathbb{R}_+)$
we have
\begin{align*}
    <\rho,\varphi_t> &+ <\alpha \rho u, \varphi_x > =
    \int_0^\infty \int_\mathbb{R} (\rho \varphi_t + \alpha(t) \rho u \varphi_x ) dx dt \\
    &+ \int_0^\infty \frac{\rho_-+\rho_+}{2}(x_-(t)-x_+(t)) \left( \varphi_t + \alpha(t) \frac{u_-+u_+}{2} \exp(-\int_0^t \sigma(s)ds) \varphi_x \right) dt \\
    =& \int_0^\infty \int_{-\infty}^{x(t)} (\rho_- \varphi_t + \alpha(t) \rho_- u_-\exp(-\int_0^t \sigma(s)ds) \varphi_x ) dx dt
+
\int_0^\infty \int_{x(t)}^\infty (\rho_+ \varphi_t + \alpha(t) \rho_+ u_+\exp(-\int_0^t \sigma(s)ds) \varphi_x ) dx dt \\
&+ \int_0^\infty \frac{\rho_-+\rho_+}{2}(x_-(t)-x_+(t)) \left( \varphi_t + \alpha(t) \frac{u_-+u_+}{2} \exp(-\int_0^t \sigma(s)ds) \varphi_x \right) dt \\
=& -\oint -(\alpha(t) \rho_- u_- \exp(-\int_0^t \sigma(s)ds) \varphi) dt + (\rho_- \varphi) dx
+ \oint -(\alpha(t) \rho_+ u_+ \exp(-\int_0^t \sigma(s)ds) \varphi) dt + (\rho_+ \varphi) dx \\
&+ \int_0^\infty \frac{\rho_-+\rho_+}{2}(x_-(t)-x_+(t)) \left( \varphi_t + \alpha(t) \frac{u_-+u_+}{2} \exp(-\int_0^t \sigma(s)ds) \varphi_x \right) dt \\
= & \int_0^t \left( \alpha(t)(\rho_- u_- - \rho_+ u_+)\exp(-\int_0^t \sigma(s)ds) - (\rho_- - \rho_+) \frac{dx(t)}{dt} \right) \varphi dt \\
&+ \int_0^\infty \frac{\rho_-+\rho_+}{2}(x_-(t)-x_+(t)) \frac{d\varphi}{dt} dt \\
= & \int_0^t \left( \alpha(t)(\rho_- u_- - \rho_+ u_+)\exp(-\int_0^t \sigma(s)ds) - (\rho_- - \rho_+) \frac{dx(t)}{dt} \right) \varphi dt \\
&- \int_0^\infty \frac{d}{dt} \left( \frac{\rho_-+\rho_+}{2}(x_-(t)-x_+(t)) \right) \varphi dt= 0, 
\end{align*}
and
\begin{align*}
    \int_0^\infty \int_\mathbb{R} & \left( u \varphi_t + \frac{\alpha(t)}{2} u^2 \varphi_x - \sigma(t) u \varphi \right) dx dt 
    = \int_0^\infty \int_\mathbb{R} \left( u \varphi_t + \frac{\alpha(t)}{2} u^2 \varphi_x \right) dx dt - \int_0^\infty \int_\mathbb{R} \sigma(t) u \varphi  dx dt \\
    =&
    \int_0^\infty \int_{-\infty}^{x(t)} u_- \exp(-\int_0^t \sigma(s)ds) \left( \varphi_t + \frac{\alpha(t)}{2} u_- \exp(-\int_0^t \sigma(s)ds) \varphi_x \right) dx dt \\
&+
\int_0^\infty \int_{x(t)}^{\infty} u_+ \exp(-\int_0^t \sigma(s)ds) \left( \varphi_t + \frac{\alpha(t)}{2} u_+ \exp(-\int_0^t \sigma(s)ds) \varphi_x \right) dx dt
- \int_0^\infty \int_\mathbb{R} \sigma(t) u \varphi  dx dt\\
=& - \oint - \left( \frac{\alpha(t)}{2}u_-^2 \exp(-2\int_0^t \sigma(s)ds) \varphi \right) dt 
+
\left( u_- \exp(-\int_0^t \sigma(s)ds ) \varphi \right) dx \\
&+ \int_0^\infty \int_{-\infty}^{x(t)} \sigma(t) u_- \exp(-\int_0^t \sigma(s)ds) \varphi dx dt \\
&+ \oint -\left( \frac{\alpha(t)}{2}u_+^2 \exp(-2\int_0^t \sigma(s)ds) \varphi \right) dt 
+
\left( u_+ \exp(-\int_0^t \sigma(s)ds ) \varphi \right) dx \\
&+ \int_0^\infty \int_{x(t)}^{\infty} \sigma(t) u_+ \exp(-\int_0^t \sigma(s)ds) \varphi dx dt - \int_0^\infty \int_\mathbb{R} \sigma(t) u \varphi  dx dt \\
=& \int_0^\infty \left( \frac{\alpha(t)}{2}(u_-^2 - u_+^2) \exp(-\int_0^t \sigma(s)ds) - (u_- - u_+) \frac{dx(t)}{dt} \right) \varphi \exp(-\int_0^t \sigma(s) ds) dt
=0.
\end{align*}

\medskip
3. Finally, we observe that, for each $t \ge 0$, 
$$
u_+ \alpha(t) \exp(-\int_0^t \sigma(\tau) d\tau ) < \frac{dx(t)}{dt} < u_+ \alpha(t) \exp(-\int_0^t \sigma(\tau) d\tau),
$$
which is an entropy condition to the system \eqref{systemLD}.
\end{proof}

\section{Pressureless Type Gas Dynamics System} 
\label{Section3}

The main issue of this section is to study the Riemann problem of the pressureless gas system with variable coefficient 
and time-variable linear damping \eqref{ZP}. We introduce a similar variable to reduce the system \eqref{ZP} to 
hyperbolic conservation laws with variable coefficient to solve the Riemann problem with $u_- < u_+$. 
To the case $u_->u_+$, similar to \cite{DelacruzJuajibioyAAM2022}, we use a nonlinear viscous system and 
using a similar variable we obtain viscous solutions that converge to a delta shock solution of the Riemann problem \eqref{ZP}-\eqref{DtIni}.

\subsection{Classical Riemann solutions.} 
\label{Sect3.1}
We observe that under transformation $\widehat{u}(x,t)=u(x,t)e^{\int_0^t \sigma(r)dr}$ the system \eqref{ZP} is equivalent to
\begin{equation} \label{ZPsysVarCoeHom}
    \begin{cases}
        \rho_t + \alpha(t) e^{-\int_0^t \sigma(r)dr} (\rho \widehat{u})_x =0,\\
        (\rho \widehat{u})_t + \alpha(t) e^{-\int_0^t \sigma(r)dr} (\rho \widehat{u}^2)_x =0,
    \end{cases}
\end{equation}
with the initial data \eqref{DtIni}.
Using the similar variable
\begin{equation} \label{similarVar}
    \xi = \frac{x}{\int_0^t \alpha(s) e^{-\int_0^s \sigma(r)dr}ds},
\end{equation}
the system \eqref{ZPsysVarCoeHom} can be written as
\begin{equation} \label{XiZPH}
    \begin{cases}
        -\xi \rho_\xi + (\rho \widehat{u})_\xi =0,\\
        -\xi (\rho \widehat{u})_\xi + (\rho \widehat{u}^2)_\xi =0,
    \end{cases}
\end{equation}
and the initial condition \eqref{DtIni} changes to the boundary condition
\begin{equation*}
    (\rho(\pm \infty),\widehat{u}(\pm \infty)) = (\rho_\pm,u_\pm).
\end{equation*}

Now, we note that any smooth solution of the system \eqref{XiZPH} satisfies
\begin{equation*}
    \begin{pmatrix}
        \widehat{u}-\xi & \rho \\
        \widehat{u}(\widehat{u}-\xi) & \rho(2\widehat{u}-\xi) 
    \end{pmatrix}
    \,
    \begin{pmatrix}
    \rho_\xi \\ \widehat{u}_\xi
    \end{pmatrix}
    =
    \begin{pmatrix}
    0 \\ 0
    \end{pmatrix}
\end{equation*}
and it provides either the general solution (constant state)
$\rho(\xi)=$constant and $\widehat{u}(\xi)=$constant, $\rho \neq 0$, or the singular solution $\rho(\xi) =0$ for all $\xi$ and $\widehat{u}(\xi)=\xi$, called the {\em vacuum state}.
Thus the smooth solutions of system \eqref{XiZPH} only contain constants and vacuum solutions. 
For a bounded discontinuity at $\xi = \eta$, the Rankine-Hugoniot condition holds, that is to say, 
\begin{equation*}
    \begin{cases}
        -\eta (\rho_--\rho_+) + (\rho_- \widehat{u}_- - \rho_+ \widehat{u}_+) = 0,\\
        -\eta (\rho_- \widehat{u}_- - \rho_+ \widehat{u}_+) + (\rho_- \widehat{u}_-^2 - \rho_+ \widehat{u}_+^2) =0,
    \end{cases}
\end{equation*}
which holds when $\eta = u_- = u_+$. Therefore, two states $(\rho_-,u_-)$ and $(\rho_+,u_+)$ can be connected by a 
contact discontinuity if and only if $u_- = u_+$. Thus, the contact discontinuity is characterized by $\xi = u_- = u_+$.

\smallskip
Summarizing, we obtain the solution which consists of two contact discontinuities and a vacuum state besides two constant states.
Therefore, the solution can be expressed as
\begin{equation*}
    (\rho(\xi),\widehat{u}(\xi))=
    \begin{cases}
        (\rho_-,u_-), &\mbox{if } \xi < u_-,\\
        (0, \xi), &\mbox{if } u_- \le \xi \le u_+,\\
        (\rho_+,u_+), &\mbox{if } \xi > u_+.
    \end{cases}
\end{equation*}
Since $u(x,t) = \widehat{u}(x,t)e^{-\int_0^t \sigma(r)dr}$ and $\xi = \frac{x}{\int_0^t \alpha(s) e^{-\int_0^s \sigma(r)dr}ds}$, then for $u_- < u_+$ the Riemann solution to the system \eqref{ZP} is
\begin{equation*}
    (\rho(x,t),u(x,t))=
    \left \{
    \begin{aligned}
        (\rho_-,u_- e^{-\int_0^t \sigma(r)dr}), &\quad \mbox{if } x < u_- \int_0^t \alpha(s)e^{-\int_0^s \sigma(r)dr}ds,\\
        (0, \frac{x e^{-\int_0^t \sigma(r)dr}}{\int_0^t \alpha(s) e^{-\int_0^s \sigma(r)dr}ds}), &\quad \mbox{if } u_- \int_0^t \alpha(s)e^{-\int_0^s \sigma(r)dr}ds \le x \le u_+ \int_0^t \alpha(s)e^{-\int_0^s \sigma(r)dr}ds,\\
        (\rho_+,u_+ e^{-\int_0^t \sigma(r)dr}), &\quad \mbox{if } x > u_+ \int_0^t \alpha(s)e^{-\int_0^s \sigma(r)dr}ds.
    \end{aligned}
    \right.
\end{equation*}

\subsection{Delta shock wave solutions.}

Given $\varepsilon> 0$, we consider the following parabolic regularization to the system \eqref{ZP},
\begin{equation} 
\label{ZeroPressureSystemVis}
\left \{
 \begin{aligned}
     \rho_t^\varepsilon + \alpha(t) (\rho^\varepsilon u^\varepsilon)_x&= 0,
     \\[5pt]
     (\rho^\varepsilon u^\varepsilon)_t + \alpha(t)(\rho^\varepsilon (u^\varepsilon)^2)_x &= \varepsilon \beta_*(t)u_{xx}^\varepsilon - \sigma(t)\rho^\varepsilon u^\varepsilon, 
 \end{aligned}
 \right.
\end{equation}
where $\beta_*(t) = \alpha(t) \exp (-\int_0^t \sigma(s) ds) \int_0^t \alpha(s) \exp (-\int_0^s \sigma(r)dr ) ds$, 
with initial condition 
\begin{equation} \label{DtVisN}
    (\rho^\varepsilon(x,0),u^\varepsilon(x,0)) = (\rho_0(x),u_0(x)),
\end{equation}
where $(\rho_0,u_0)$ is given by \eqref{DtIni}. 

\medskip
Under the transformation
$\widehat{u}^\varepsilon(x,t)=u^\varepsilon(x,t)e^{\int_0^t \sigma(r)dr}$ the system \eqref{ZeroPressureSystemVis} becomes 
\begin{equation} 
\label{systemVis4}
\left \{
    \begin{aligned}
        \rho_t^\varepsilon + \alpha(t)e^{-\int_0^t \sigma(r)dr} (\rho^\varepsilon \widehat{u}^\varepsilon)_x&=0,
        \\[5pt]
        (\rho^\varepsilon \widehat{u}^\varepsilon)_t + \alpha(t) e^{-\int_0^t \sigma(r)dr} (\rho^\varepsilon (\widehat{u}^\varepsilon)^2)_x&= \varepsilon \beta_*(t) \widehat{u}_{xx}^\varepsilon, 
    \end{aligned}
    \right. 
\end{equation}
and the initial condition \eqref{DtVisN} becomes
\begin{equation} \label{DtVisN2}
    (\rho^\varepsilon(x,0),\widehat{u}^\varepsilon(x,0)) = (\rho_0^\varepsilon(x),\widehat{u}_0^\varepsilon(x)) =
    \begin{cases}
    (\rho_-,u_-), &\mbox{if } x<0,\\
    (\rho_+,u_+), &\mbox{if } x>0
    \end{cases}
\end{equation}
for arbitrary constant states $u_\pm$ and $\rho_\pm >0$ as well. 
By using the similar variable \eqref{similarVar}
the system \eqref{systemVis4} 
can be written as
\begin{equation} 
\label{systemVis5}
\left \{
    \begin{aligned}
        -\xi \rho_\xi^\varepsilon + (\rho^\varepsilon \widehat{u}^\varepsilon)_\xi&=0,
        \\[5pt]
        - \xi (\rho^\varepsilon \widehat{u}^\varepsilon)_\xi + (\rho^\varepsilon (\widehat{u}^\varepsilon)^2)_{\xi}&= \varepsilon \widehat{u}_{\xi \xi}^\varepsilon
    \end{aligned}
    \right. 
\end{equation}
and the initial data \eqref{DtVisN2} changes to the boundary condition
\begin{equation} \label{BoundCondZP}
(\rho(\pm \infty), \widehat{u}(\pm \infty))=(\rho_\pm,u_\pm)
\end{equation}
for arbitrary constant states $u_->u_+$ and $\rho_\pm >0$.
The existence of solutions to the system \eqref{systemVis5} with boundary condition \eqref{BoundCondZP} was shown in Theorem 3 of \cite{DelacruzJuajibioyAAM2022}. More explicitly,  in \cite{DelacruzJuajibioyAAM2022}, the following result was obtained:
\begin{proposition}
    There exists a weak solution $(\rho^\varepsilon, \widehat{u}^\varepsilon) \in L_{loc}^1((-\infty,+\infty)) \times C^2((-\infty,+\infty))$ to the boundary problem \eqref{systemVis5}-\eqref{BoundCondZP}.
\end{proposition}

From Theorem 2 in \cite{DelacruzJuajibioyAAM2022}, we have that for each $\varepsilon >0$, the function $\widehat{u}^\varepsilon$ satisfies
\begin{equation*}
    \begin{cases}
        \varepsilon (\widehat{u}^\varepsilon)''(\xi) = (\rho^\varepsilon(\xi)(\widehat{u}-\xi)) (\widehat{u}^\varepsilon)'(\xi),\\
        \widehat{u}^\varepsilon(\pm \infty) = u_\pm,
    \end{cases}
\end{equation*}
with $' = \frac{d}{d\xi}$ and 
\begin{equation*}
    \rho^\varepsilon(\xi) =
    \begin{cases}
        \rho_1^\varepsilon(\xi), &\mbox{if } -\infty<\xi <\xi_\varsigma^\varepsilon,\\
        \rho_2^\varepsilon(\xi), &\mbox{if } \xi_\varsigma^\varepsilon < \xi < +\infty,
    \end{cases}
\end{equation*}
where $\xi_\varsigma^\varepsilon$ satisfies $\widehat{u}^\varepsilon(\xi_\varsigma^\varepsilon) = \xi_\varsigma^\varepsilon$,
\begin{equation*}
    \rho_1(\xi) = \rho_-\exp \left( - \int_{-\infty}^\xi \frac{(\widehat{u}^\varepsilon(s))'}{\widehat{u}^\varepsilon(s)-s}ds \right)
    \quad
    \mbox{ and }
    \quad
    \rho_2(\xi) = \rho_+\exp \left( \int_\xi^{\infty} \frac{(\widehat{u}^\varepsilon(s))'}{\widehat{u}^\varepsilon(s)-s}ds \right).
\end{equation*}

\begin{definition} 
    A distribution pair $(\rho,u)$ is called a delta shock wave solution of the problem \eqref{ZP} and \eqref{DtIni} in the sense of distributions, 
    when there exist a smooth curve $L$ and a function $w(\cdot)$, such that $\rho$ and $u$ are represented respectively by 
    $$
    \rho= \widehat{\rho}(x,t) + w\delta_L \quad \mbox{ and } \quad u= u(x,t),
    $$
with $\widehat{\rho}, u \in L^\infty(\mathbb{R} \times (0,\infty))$, and satisfy
     for each the test function $\varphi \in C_0^\infty(\mathbb{R} \times (0,\infty))$,
    \begin{equation} 
    \label{DelSolZPsys}
        \begin{cases}
            < \rho,\varphi_t>+<\alpha \rho u,\varphi_x>=0,\\
            < \rho u,\varphi_t> + <\alpha \rho u^2,\varphi_x > = < \sigma \rho u, \varphi >,
        \end{cases}
    \end{equation}
    where
    $$
    <\rho,\varphi>:= \iint_{\mathbb{R}^2_+} \widehat{\rho} \varphi dx dt + \langle w \delta_L, \varphi \rangle
    $$
    and for some smooth function $G$, 
    $$
    <\alpha \rho G(\cdot),\varphi>:= \iint_{\mathbb{R}^2_+} \alpha(t) \widehat{\rho} G(u) \varphi dx dt + \langle \alpha(\cdot) w G(u_\delta) \delta_L, \varphi \rangle. 
    $$
  Moreover, $u|_L = u_\delta(t)$. 
\end{definition}

\medskip
Now, we denote $\varsigma = \lim\limits_{\varepsilon \to 0+} \xi_\varsigma^\varepsilon = \lim\limits_{\varepsilon \to 0+} \widehat{u}^\varepsilon(\xi_\varsigma^\varepsilon) = \widehat{u}(\varsigma)$.
Then, according to Theorem 4 in \cite{DelacruzJuajibioyAAM2022}, we have
\begin{equation*} 
\lim\limits_{\varepsilon \to 0+} (\rho^\varepsilon(\xi), \widehat{u}^\varepsilon(\xi))=
    \begin{cases}
        (\rho_-,u_-), &\mbox{if } \xi < \varsigma,\\
(w_0 \, \delta(\xi-\varsigma), u_\delta), &\mbox{if } \xi = \varsigma,\\
(\rho_+,u_+), &\mbox{if } \xi > \varsigma,
    \end{cases}
\end{equation*}
where $\rho^\varepsilon$ converges in the sense of distributions to the sum of a step function and a Dirac measure $\delta$ 
with weight $w_0 = -\varsigma (\rho_--\rho_+) + (\rho_-u_- - \rho_+ u_+)$ and $u_\delta = \widehat{u}(\varsigma)$. 
Moreover, $(\varsigma,w_0,u_\delta)$ satisfies
\begin{equation} \label{RHZP1}
    \begin{cases}
        \varsigma = u_\delta,\\
        w_0 = -\varsigma(\rho_--\rho_+) + (\rho_-u_- - \rho_+u_+),\\
        w_0 u_\delta = - \varsigma (\rho_-u_- - \rho_+u_+) + (\rho_-u_-^2 - \rho_+u_+^2),
    \end{cases}
\end{equation}
and the over-compressive entropy condition 
\begin{equation} \label{entropy1}
u_+<u_\delta<u_-.
\end{equation}

Observe that from the system \eqref{RHZP1} we have
$$
(\rho_--\rho_+)u_\delta^2 - 2 (\rho_- u_- - \rho_+ u_+)u_\delta + (\rho_- u_-^2 - \rho_+ u_+^2)= 0, 
$$
which implies
$$
u_\delta = \frac{\sqrt{\rho_-}u_--\sqrt{\rho_+}u_+}{\sqrt{\rho_-} - \sqrt{\rho_+}}
\quad
\mbox{ or }
\quad
u_\delta = \frac{\sqrt{\rho_-}u_-+\sqrt{\rho_+}u_+}{\sqrt{\rho_-} + \sqrt{\rho_+}}.
$$

One remarks that, when $u_\delta = \frac{\sqrt{\rho_-}u_-+\sqrt{\rho_+}u_+}{\sqrt{\rho_-} + \sqrt{\rho_+}}$
the entropy condition is valid while $u_\delta = \frac{\sqrt{\rho_-}u_--\sqrt{\rho_+}u_+}{\sqrt{\rho_-} - \sqrt{\rho_+}}$ does not satisfy the entropy condition.
Moreover, using the second equation of the system \eqref{RHZP1} and $u_\delta = \frac{\sqrt{\rho_-}u_-+\sqrt{\rho_+}u_+}{\sqrt{\rho_-} + \sqrt{\rho_+}}$, we obtain $w_0 = \sqrt{\rho_- \rho_+}(u_--u_+)$. 
Therefore, when $\rho_-=\rho_+$, from \eqref{RHZP1} we obtain 
$$
    2(u_--u_+)u_\delta -(u_-^2 - u_+^2)= 0
$$ 
and hence we have $u_\delta= \frac{1}{2}(u_-+u_+)$ and $w_0=\rho_-(u_--u_+)$.
Finally, using the similar variable \eqref{similarVar}, we have obtained the following result
\begin{proposition}
    Suppose $u_->u_+$. Let $(\rho^\varepsilon(x,t), \widehat{u}^\varepsilon(x,t))$ be the solution of the problem \eqref{systemVis4}-\eqref{DtVisN2}.
    Then the limit $\lim\limits_{\varepsilon \to 0+} (\rho^\varepsilon(x,t), \widehat{u}^\varepsilon(x,t)) = (\rho(x,t), \widehat{u}(x,t))$ exists in the distribution sense. Moreover,
    $ (\rho(x,t), \widehat{u}(x,t))$ is given by 
    \begin{equation*}
 \left \{
        \begin{aligned}
            (\rho_-,u_-) \quad , &\quad \mbox{if } x < u_\delta \int_0^t \alpha(s) e^{-\int_0^s \sigma(r) dr} ds,
            \\
            (w_0 \int_0^t \alpha(s) e^{-\int_0^s \sigma(r) dr} ds \cdot \delta(x-u_\delta \int_0^t \alpha(s) e^{-\int_0^s \sigma(r) dr} ds), u_\delta), &\quad \mbox{if } x = u_\delta \int_0^t \alpha(s) e^{-\int_0^s \sigma(r) dr} ds,
            \\
            (\rho_+,u_+) \quad , &\quad \mbox{if } x > u_\delta \int_0^t \alpha(s) e^{-\int_0^s \sigma(r) dr} ds,
        \end{aligned}
        \right.
    \end{equation*}
    where $w_0 = \sqrt{\rho_- \rho_+}(u_--u_+)$ and $u_\delta = \frac{\sqrt{\rho_-}u_-+\sqrt{\rho_+}u_+}{\sqrt{\rho_-} + \sqrt{\rho_+}}$, when $\rho_- \neq \rho_+$. 
    For the case $\rho_- = \rho_+$, it follows that, $w_0 = \rho_-(u_--u_+)$ and $u_\delta = \frac{1}{2}(u_-+u_+)$.
    In addition, the solution is unique under the over-compressive entropy condition \eqref{entropy1}.
\end{proposition}

\begin{remark}
The condition \eqref{RHZP1} is necessary and sufficient to guarantee the existence of delta shock solutions to the problem \eqref{systemVis4}-\eqref{DtVisN2} with $\varepsilon =0$. 
In fact, there are two delta shock solutions. 
Now, the over-compressive entropy condition \eqref{entropy1}, (see the above proposition), was sufficient to obtain the uniqueness of the delta shock solution.
\end{remark}

\medskip
From the above proposition and since $u(x,t) = \widehat{u}(x,t) e^{-\int_0^t \sigma(r)dr}$, we can establish 
a solution to the system \eqref{ZP} with initial data \eqref{DtIni}. 
Moreover, multiplying the entropy condition \eqref{entropy1} by $\alpha(t)$ we get $\alpha(t)u_+ < u_\delta \alpha(t) < \alpha(t)$ for all $t \ge 0$ and again using that
$u(x,t) = \widehat{u}(x,t) e^{-\int_0^t \sigma(r)dr}$,
we have extended the entropy condition \eqref{entropy1} to the following entropy condition to the system \eqref{ZP},
\begin{equation} \label{entropy2}
\lambda(\rho_+,u_+) e^{-\int_0^t \sigma(r)dr}  < \frac{dx(t)}{dt} < \lambda(\rho_-,u_-)e^{-\int_0^t \sigma(r)dr}, \mbox{ for all } t \ge 0,
\end{equation}
where $\lambda(\rho,u) = \alpha u$ is the eigenvalue associated to system \eqref{ZP}.
Then, we have the following
\begin{theorem}
    Suppose $u_->u_+$. Then the Riemann problem \eqref{ZP}-\eqref{DtIni} admits 
    under the entropy condition \eqref{entropy2} a unique delta shock solution of the form
    \begin{equation} 
    \label{DeltZPs}
        (\rho(x,t),u(x,t))=
        \left \{
        \begin{aligned}
            (\rho_-,u_- e^{-\int_0^t \sigma(r)dr}), &\quad \mbox{if } x < x(t),\\
            (w(t) \delta(x-x(t)), u_\delta(t)), &\quad \mbox{if } x = x(t),\\
            (\rho_+,u_+ e^{-\int_0^t \sigma(r)dr}), &\quad \mbox{if } x > x(t),
        \end{aligned}
        \right.
    \end{equation}
    where for $\rho_- \neq \rho_+$, 
    $$
    \begin{aligned}
        w(t)&= \sqrt{\rho_- \rho_+} \, (u_--u_+) \int_0^t \alpha(s) e^{-\int_0^s \sigma(r) dr} ds, \quad 
        u_\delta(t)= \frac{\sqrt{\rho_-}u_-+\sqrt{\rho_+}u_+}{\sqrt{\rho_-} + \sqrt{\rho_+}} e^{-\int_0^t \sigma(r)dr}, \quad \text{and} 
        \\[5pt] 
        x(t)&= \frac{\sqrt{\rho_-}u_-+\sqrt{\rho_+}u_+}{\sqrt{\rho_-} + \sqrt{\rho_+}} \int_0^t \alpha(s) e^{-\int_0^s \sigma(r) dr} ds. 
\end{aligned}  
$$      
For the case $\rho_- = \rho_+$, it follws that 
$$
\begin{aligned}
w(t)&=\rho_-(u_--u_+) \int_0^t \alpha(s) e^{-\int_0^s \sigma(r) dr} ds, \quad u_\delta(t) = \frac{1}{2}(u_-+u_+) e^{-\int_0^s \sigma(r) dr}, \quad  \text{and} 
\\[5pt]
x(t)&= \frac{1}{2}(u_-+u_+) \int_0^t \alpha(s) e^{-\int_0^s \sigma(r) dr} ds.
\end{aligned}
$$
\end{theorem}
\begin{proof}
Suppose that $\rho_- \neq \rho_+$. Therefore, in order to show that $(\rho,u)$, given by \eqref{DeltZPs}, 
is a solution to the problem \eqref{ZP}-\eqref{DtIni}, we consider any test function $\varphi \in C_0^\infty( \mathbb{R} \times (0,\infty))$ and compute 
    $$
    \begin{aligned}
        <\rho u, \varphi_t> &+ <\rho u^2, \varphi_x >  
        = \int_0^\infty \int_\mathbb{R} (\rho u \varphi_t + \alpha(t) \rho u^2 \varphi_x) dx dt
        + \int_0^\infty  w(t) (u_\delta(t) \varphi_t + \alpha(t) u^2_\delta(t) \varphi_x) dt \\
        =& \int_0^\infty \int_{-\infty}^{x(t)} (\rho_- u_- e^{-\int_0^t \sigma(r)dr} \varphi_t + \alpha(t) \rho_- u_-^2 e^{-2\int_0^t \sigma(r)dr} \varphi_x) dx dt \\
        &+ \int_0^\infty \int_{x(t)}^{\infty} (\rho_+ u_+ e^{-\int_0^t \sigma(r)dr} \varphi_t + \alpha(t) \rho_+ u_+^2 e^{-2\int_0^t \sigma(r)dr} \varphi_x) dx dt\\
        &+ \int_0^\infty w(t) \frac{\sqrt{\rho_-}u_-+\sqrt{\rho_+}u_+}{\sqrt{\rho_-} + \sqrt{\rho_+}} e^{-\int_0^t \sigma(r)dr} \left( \varphi_t + \alpha(t) \frac{\sqrt{\rho_-}u_-+\sqrt{\rho_+}u_+}{\sqrt{\rho_-} + \sqrt{\rho_+}} e^{-\int_0^t \sigma(r)dr} \varphi_x \right) dt \\
        =& - \oint - \left( \alpha(t) \rho_- u_-^2 e^{-2\int_0^t \sigma(r)dr} \varphi \right) dt + \left( \rho_- u_- e^{-\int_0^t \sigma(r)dr} \varphi \right) dx \\
        &+ \oint - \left( \alpha(t) \rho_+ u_+^2 e^{-2\int_0^t \sigma(r)dr} \varphi \right) dt + \left( \rho_+ u_+ e^{-\int_0^t \sigma(r)dr} \varphi \right) dx\\
        &+ \int_0^\infty \int_\mathbb{R} \sigma(t) \rho u \varphi dx dt 
        + \int_0^\infty w(t) \frac{\sqrt{\rho_-}u_-+\sqrt{\rho_+}u_+}{\sqrt{\rho_-} + \sqrt{\rho_+}} e^{-\int_0^t \sigma(r)dr} \left( \varphi_t + \frac{dx(t)}{dt} \varphi_x \right) dt 
                    \end{aligned}
    $$
        $$
    \begin{aligned}
        =& \int_0^\infty \alpha(t) (\rho_- u_-^2 - \rho_+ u_+^2) e^{-2\int_0^t \sigma(r)dr} \varphi dt
        - \int_0^\infty \frac{dx(t)}{dt} (\rho_- u_- - \rho_+ u_+) e^{-\int_0^t \sigma(r)dr} \varphi dt \\
        & + \int_0^\infty \int_\mathbb{R} \sigma(t) \rho u \varphi dx dt 
        + \int_0^\infty w(t) \frac{\sqrt{\rho_-}u_-+\sqrt{\rho_+}u_+}{\sqrt{\rho_-} + \sqrt{\rho_+}} e^{-\int_0^t \sigma(r)dr} \frac{d\varphi(t)}{dt} dt \\
        =& \int_0^\infty \alpha(t) (\rho_- u_-^2 - \rho_+ u_+^2) e^{-2\int_0^t \sigma(r)dr} \varphi dt
        - \int_0^\infty \frac{dx(t)}{dt} (\rho_- u_- - \rho_+ u_+) e^{-\int_0^t \sigma(r)dr} \varphi dt \\
        & + \int_0^\infty \int_\mathbb{R} \sigma(t) \rho u \varphi dx dt 
        - \int_0^\infty \frac{d}{dt} \left( w(t) \frac{\sqrt{\rho_-}u_-+\sqrt{\rho_+}u_+}{\sqrt{\rho_-} + \sqrt{\rho_+}} e^{-\int_0^t \sigma(r)dr} \right) \varphi dt \\
        =& \int_0^\infty \int_\mathbb{R} \sigma(t) \rho u \varphi dx dt 
        + \int_0^\infty \sigma(t) w(t) \frac{\sqrt{\rho_-}u_-+\sqrt{\rho_+}u_+}{\sqrt{\rho_-} + \sqrt{\rho_+}} e^{-\int_0^t \sigma(r)dr} dt
        =  < \sigma \rho u, \varphi >, 
    \end{aligned}
    $$
    which implies the second equation of \eqref{DelSolZPsys}.
    With a similar argument, it is possible to obtain the first equation of \eqref{DelSolZPsys} and the case when $\rho_- = \rho_+$.
    The uniqueness of the solution will be obtained under the entropy condition \eqref{entropy2}. 
\end{proof}

\section{Riemann problem to the systems (4) and (5) with $\sigma(\cdot) \equiv 0$}

In this section, we consider $\sigma(t) = \mu \nu(t)$ where $\mu >0$ is a parameter, $\nu(t)\geq 0$ for all $t \ge 0$, and $\nu \in L_{\rm loc}^1([0,\infty))$.
According to the Sections \ref{Sect2.3.1} and \ref{Sect3.1}, if $u_- < u_+$, the systems \eqref{systemLD} and \eqref{ZP}
with initial data \eqref{DtIni} have the solution
\begin{equation*}
    (\rho(x,t),u(x,t)) =
    \begin{cases}
        (\rho_-, u_- \exp(-\mu \int_0^t \nu(s)ds)), &\mbox{if } x < x_-(t),\\
        (0, \frac{x}{\int_0^t \alpha(s) \exp(-\mu \int_0^s \nu(r)dr)ds} \exp(-\mu \int_0^t \nu(s)ds)), &\mbox{if } x_-(t) <x< x_+(t),\\
        (\rho_+, u_+ \exp(-\mu \int_0^t \nu(s)ds)), &\mbox{if } x > x_+(t),\\
    \end{cases}
\end{equation*}
where $x_\pm(t) = u_\pm \int_0^t \alpha(s) \exp(-\mu \int_0^s \nu(r)dr) ds$. 
If $u_->u_+$, the the solution for the problem \eqref{systemLD}-\eqref{DtIni} is 
\begin{equation*}
    (\rho(x,t),u(x,t)) =
    \begin{cases}
        (\rho_-,u_- \exp (- \mu \int_0^t \nu(s)ds)), &\mbox{if } x<x(t),\\
        (w(t)\delta(x-x(t)), \frac{u_-+u_+}{2} \exp (- \mu \int_0^t \nu(s)ds)), &\mbox{if } x=x(t),\\
        (\rho_+,u_+ \exp (- \mu \int_0^t \nu(s)ds)), &\mbox{if } x>x(t),
    \end{cases}
\end{equation*}
where $w(t) = \frac{1}{2}(\rho_-+\rho_+)(u_--u_+) \int_0^t \alpha(s) \exp (- \mu \int_0^s \nu(\tau) d\tau) ds$ and $x(t) = \frac{u_-+u_+}{2} \int_0^t \alpha(s) \exp (- \mu \int_0^s \nu(\tau) d\tau) ds$ while the solution to the problem \eqref{ZP}-\eqref{DtIni}
is
\begin{equation*}
    (\rho(x,t),u(x,t))=
    \begin{cases}
        (\rho_-,u_-e^{-\mu \int_0^t \nu(r)dr}, &\mbox{if } x < x(t),\\
        (w(t)\delta(x-x(t)),u_\delta(t)), &\mbox{if } x=x(t),\\
        (\rho_+,u_+e^{-\mu \int_0^t \nu(r)dr}, &\mbox{if } x > x(t),
    \end{cases}
\end{equation*}
where $w(t) = \sqrt{\rho_+ \rho_-}(u_- - u_+) \int_0^t \alpha(s) e^{-\mu \int_0^s \nu(r)dr}ds$, $u_\delta(t) = \frac{\sqrt{\rho_-}u_- + \sqrt{\rho_+}u_+}{\sqrt{\rho_-} + \sqrt{\rho_+}}e^{-\mu \int_0^t \nu(r)dr}$,\\ and $x(t)= \frac{\sqrt{\rho_-}u_- + \sqrt{\rho_+}u_+}{\sqrt{\rho_-} + \sqrt{\rho_+}} \int_0^t \alpha(s) e^{-\mu \int_0^s \nu(r)dr}ds$ if $\rho_- \neq \rho_+$ and $w(t) = \rho_-(u_- - u_+) \int_0^t \alpha(s) e^{-\mu \int_0^s \nu(r)dr}ds$, $u_\delta(t) = \frac{1}{2}(u_- + u_+) e^{-\mu \int_0^t \nu(r)dr}$, and $x(t)=\frac{1}{2}(u_- - u_+) \int_0^t \alpha(s) e^{-\mu \int_0^s \nu(r)dr}ds$ if $\rho_- = \rho_+$.\\

One observes that the solutions given above are explicit with respect to the parameter $\mu> 0$, and also we have
$$ 
     \lim\limits_{\mu \to 0+} \exp(-\mu \int_0^t \nu(s)ds) = 1
\quad \text{and} \quad 
\lim\limits_{\mu \to 0+} \int_0^t \alpha(s) \exp(-\mu \int_0^s \nu(r)dr) ds = \int_0^t \alpha(s) ds.
$$
Therefore, the Riemann solution to the problems \eqref{systemLD} and \eqref{ZP} with $\sigma(t)=0$ for all $t \ge 0$ and initial data \eqref{DtIni} is given by
\begin{equation*}
    (\rho(x,t),u(x,t)) =
    \begin{cases}
        (\rho_-, u_-), &\mbox{if } x < u_- \int_0^t \alpha(s)ds,\\
        (0, \frac{x}{\int_0^t \alpha(s) ds} ), &\mbox{if } u_- \int_0^t \alpha(s)ds <x< u_+ \int_0^t \alpha(s)ds,\\
        (\rho_+, u_+), &\mbox{if } x > u_+ \int_0^t \alpha(s)ds.
    \end{cases}
\end{equation*}
if $u_- < u_+$. If $u_->u_+$, then the Riemann solution to the problem \eqref{systemLD} with $\sigma(t)=0$ for all $t \ge 0$ and initial data \eqref{DtIni} is
\begin{equation*}
    (\rho(x,t),u(x,t)) =
    \begin{cases}
        (\rho_-,u_-), &\mbox{if } x<x(t),\\
        (w(t)\delta(x-x(t)), \frac{u_-+u_+}{2}), &\mbox{if } x=x(t),
        \\
        (\rho_+,u_+), &\mbox{if } x>x(t),
    \end{cases}
\end{equation*}
where $w(t) = \frac{1}{2}(\rho_-+\rho_+)(u_--u_+) \int_0^t \alpha(s) ds$ and $x(t) = \frac{u_-+u_+}{2} \int_0^t \alpha(s) ds$ and
the Riemann solution to the problem \eqref{ZP}-\eqref{DtIni} with $\sigma(t)=0$ for all $t \ge 0$ is given by
\begin{equation*}
    (\rho(x,t),u(x,t)) =
    \begin{cases}
        (\rho_-,u_-), &\mbox{if } x<x(t),\\
        (w(t)\delta(x-x(t)),u_\delta(t)), &\mbox{if } x=x(t),\\
        (\rho_+,u_+), &\mbox{if } x>x(t),
    \end{cases}
\end{equation*}
where $w(t)=\sqrt{\rho_+ \rho_-} (u_--u_+) \int_0^t \alpha(s)ds$, $u_\delta(t)= \frac{\sqrt{\rho_-}u_- + \sqrt{\rho_+}u_+}{\sqrt{\rho_-} + \sqrt{\rho_+}}$, and $x(t)= \frac{\sqrt{\rho_-}u_- + \sqrt{\rho_+}u_+}{\sqrt{\rho_-} + \sqrt{\rho_+}} \int_0^t \alpha(s) ds$ if $\rho_- \neq \rho_+$ and $w(t) = \rho_-(u_- - u_+) \int_0^t \alpha(s) ds$, $u_\delta(t) = \frac{1}{2}(u_- + u_+)$, and $x(t)=\frac{1}{2}(u_- - u_+) \int_0^t \alpha(s) ds$ if $\rho_- = \rho_+$.

\section{Comments and Extensions}
\label{Comments}

The main goal of this section is to present comments and extensions of ongoing work on the topic developed in this paper. 


We studied in this paper, the Riemann problems to the time-variable coefficient Zeldovich approximate system \eqref{systemLD} 
and time-variable coefficient pressureless gas system \eqref{ZP} both with general time-gradually- degenerate damping.
Similar to the results obtained by Keita and Bourgault in \cite{Keita} 
to the Riemann problems \eqref{systemLD}-\eqref{DtIni} and \eqref{ZP}-\eqref{DtIni} both with $\alpha(\cdot) \equiv 1$ 
and $\sigma(\cdot) \equiv \sigma = const.>0$, we have that the systems \eqref{systemLD} and \eqref{ZP}, where $\alpha$ 
and $\sigma$ are non-negative functions that dependents of time $t$, are equivalent for smooth and two-contact-discontinuity 
solutions but they differ for delta shock solutions.
Moreover, we show that the uniqueness is obtained under an over-compressive entropy condition.

\medskip
It is interesting to remark that,  why we have to fix the sign of $\alpha(\cdot)$ solving the Riemann problem. Indeed, they only need to have one sign (positive or negative) 
to maintain the Lax entropy (in shocks) and over-compressive entropy condition in delta shocks (as we need the characteristics not to be inverted). Clearly, the sign
of $\sigma(\cdot)$ justifies the physical meaning of damping. 

\bigskip
Now, we would like to mention another direction of the work developed here, see \cite{Richard_Wladimir}. 
Also related to system \eqref{EQQF}, we 
consider the following nonautonomous quasilinear systems:
\begin{equation*} 
\label{systemLDV}
\left \{
\begin{aligned}
    \rho_t + \alpha_1(t) \, (\rho \, u)_x&= 0,
    \\
    u_t + \alpha_2(t) \, (\frac{u^2}{2})_x&= - \sigma(t) \, u,
\end{aligned}
\right.
\end{equation*}
and also 
\begin{equation*} 
\label{ZPV}
\left \{
\begin{aligned}
    \rho_t + \alpha_1(t) \, (\rho \, u)_x&= 0,
    \\
    (\rho u)_t + \alpha_2(t) \, (\rho u^2)_x&= - \sigma(t) \, \rho u,
\end{aligned}
\right.
\end{equation*}
where $\alpha_i \in L^1([0,\infty))$, $(i= 1,2)$, and $0 \leq \sigma \in L_{\rm loc}^1([0,\infty))$. 
It is not absolutely clear that, all the strategies applied in this paper work with these systems,
in fact, this is not the case. Indeed, when $\alpha_1 \neq \alpha_2$ the construction of shocks, 
rarefactions, contact discontinuities, and delta shock solutions is not easy due to the behavior 
of the under- or over-compressibility of the eigenvalues and left or right states. 
This stands as the focal point of our ongoing research efforts. 

\section*{Data availability statement}

Data sharing does not apply to this article as no data sets were generated or analyzed during the current study.

 \section*{Conflict of Interest}

The author Richard De la cruz acknowledges the support received from Universidad Pedag{\'o}gica y Tecnol{\'o}gica de Colombia. 
The author Wladimir Neves has received research grants from CNPq
through the grants  313005/2023-0, 406460/2023-0, and also by FAPERJ 
(Cientista do Nosso Estado) through the grant E-26/201.139/2021.

\end{document}